
\catcode`\@=11


\message{Loading jyTeX fonts...}



\font\vptrm=cmr5 \font\vptmit=cmmi5 \font\vptsy=cmsy5 \font\vptbf=cmbx5

\skewchar\vptmit='177 \skewchar\vptsy='60 \fontdimen16
\vptsy=\the\fontdimen17 \vptsy

\def\vpt{\ifmmode\err@badsizechange\else
     \@mathfontinit
     \textfont0=\vptrm  \scriptfont0=\vptrm  \scriptscriptfont0=\vptrm
     \textfont1=\vptmit \scriptfont1=\vptmit \scriptscriptfont1=\vptmit
     \textfont2=\vptsy  \scriptfont2=\vptsy  \scriptscriptfont2=\vptsy
     \textfont3=\xptex  \scriptfont3=\xptex  \scriptscriptfont3=\xptex
     \textfont\bffam=\vptbf
     \scriptfont\bffam=\vptbf
     \scriptscriptfont\bffam=\vptbf
     \@fontstyleinit
     \def\rm{\vptrm\fam=\z@}%
     \def\bf{\vptbf\fam=\bffam}%
     \def\oldstyle{\vptmit\fam=\@ne}%
     \rm\fi}


\font\viptrm=cmr6 \font\viptmit=cmmi6 \font\viptsy=cmsy6
\font\viptbf=cmbx6

\skewchar\viptmit='177 \skewchar\viptsy='60 \fontdimen16
\viptsy=\the\fontdimen17 \viptsy

\def\vipt{\ifmmode\err@badsizechange\else
     \@mathfontinit
     \textfont0=\viptrm  \scriptfont0=\vptrm  \scriptscriptfont0=\vptrm
     \textfont1=\viptmit \scriptfont1=\vptmit \scriptscriptfont1=\vptmit
     \textfont2=\viptsy  \scriptfont2=\vptsy  \scriptscriptfont2=\vptsy
     \textfont3=\xptex   \scriptfont3=\xptex  \scriptscriptfont3=\xptex
     \textfont\bffam=\viptbf
     \scriptfont\bffam=\vptbf
     \scriptscriptfont\bffam=\vptbf
     \@fontstyleinit
     \def\rm{\viptrm\fam=\z@}%
     \def\bf{\viptbf\fam=\bffam}%
     \def\oldstyle{\viptmit\fam=\@ne}%
     \rm\fi}

\font\viiptrm=cmr7 \font\viiptmit=cmmi7 \font\viiptsy=cmsy7
\font\viiptit=cmti7 \font\viiptbf=cmbx7

\skewchar\viiptmit='177 \skewchar\viiptsy='60 \fontdimen16
\viiptsy=\the\fontdimen17 \viiptsy

\def\viipt{\ifmmode\err@badsizechange\else
     \@mathfontinit
     \textfont0=\viiptrm  \scriptfont0=\vptrm  \scriptscriptfont0=\vptrm
     \textfont1=\viiptmit \scriptfont1=\vptmit \scriptscriptfont1=\vptmit
     \textfont2=\viiptsy  \scriptfont2=\vptsy  \scriptscriptfont2=\vptsy
     \textfont3=\xptex    \scriptfont3=\xptex  \scriptscriptfont3=\xptex
     \textfont\itfam=\viiptit
     \scriptfont\itfam=\viiptit
     \scriptscriptfont\itfam=\viiptit
     \textfont\bffam=\viiptbf
     \scriptfont\bffam=\vptbf
     \scriptscriptfont\bffam=\vptbf
     \@fontstyleinit
     \def\rm{\viiptrm\fam=\z@}%
     \def\it{\viiptit\fam=\itfam}%
     \def\bf{\viiptbf\fam=\bffam}%
     \def\oldstyle{\viiptmit\fam=\@ne}%
     \rm\fi}


\font\viiiptrm=cmr8 \font\viiiptmit=cmmi8 \font\viiiptsy=cmsy8
\font\viiiptit=cmti8
\font\viiiptbf=cmbx8

\skewchar\viiiptmit='177 \skewchar\viiiptsy='60 \fontdimen16
\viiiptsy=\the\fontdimen17 \viiiptsy

\def\viiipt{\ifmmode\err@badsizechange\else
     \@mathfontinit
     \textfont0=\viiiptrm  \scriptfont0=\viptrm  \scriptscriptfont0=\vptrm
     \textfont1=\viiiptmit \scriptfont1=\viptmit \scriptscriptfont1=\vptmit
     \textfont2=\viiiptsy  \scriptfont2=\viptsy  \scriptscriptfont2=\vptsy
     \textfont3=\xptex     \scriptfont3=\xptex   \scriptscriptfont3=\xptex
     \textfont\itfam=\viiiptit
     \scriptfont\itfam=\viiptit
     \scriptscriptfont\itfam=\viiptit
     \textfont\bffam=\viiiptbf
     \scriptfont\bffam=\viptbf
     \scriptscriptfont\bffam=\vptbf
     \@fontstyleinit
     \def\rm{\viiiptrm\fam=\z@}%
     \def\it{\viiiptit\fam=\itfam}%
     \def\bf{\viiiptbf\fam=\bffam}%
     \def\oldstyle{\viiiptmit\fam=\@ne}%
     \rm\fi}


\def\getixpt{%
     \font\ixptrm=cmr9
     \font\ixptmit=cmmi9
     \font\ixptsy=cmsy9
     \font\ixptit=cmti9
     \font\ixptbf=cmbx9
     \skewchar\ixptmit='177 \skewchar\ixptsy='60
     \fontdimen16 \ixptsy=\the\fontdimen17 \ixptsy}

\def\ixpt{\ifmmode\err@badsizechange\else
     \@mathfontinit
     \textfont0=\ixptrm  \scriptfont0=\viiptrm  \scriptscriptfont0=\vptrm
     \textfont1=\ixptmit \scriptfont1=\viiptmit \scriptscriptfont1=\vptmit
     \textfont2=\ixptsy  \scriptfont2=\viiptsy  \scriptscriptfont2=\vptsy
     \textfont3=\xptex   \scriptfont3=\xptex    \scriptscriptfont3=\xptex
     \textfont\itfam=\ixptit
     \scriptfont\itfam=\viiptit
     \scriptscriptfont\itfam=\viiptit
     \textfont\bffam=\ixptbf
     \scriptfont\bffam=\viiptbf
     \scriptscriptfont\bffam=\vptbf
     \@fontstyleinit
     \def\rm{\ixptrm\fam=\z@}%
     \def\it{\ixptit\fam=\itfam}%
     \def\bf{\ixptbf\fam=\bffam}%
     \def\oldstyle{\ixptmit\fam=\@ne}%
     \rm\fi}


\font\xptrm=cmr10 \font\xptmit=cmmi10 \font\xptsy=cmsy10
\font\xptex=cmex10 \font\xptit=cmti10 \font\xptsl=cmsl10
\font\xptbf=cmbx10 \font\xpttt=cmtt10 \font\xptss=cmss10
\font\xptsc=cmcsc10 \font\xptbfs=cmb10 \font\xptbmit=cmmib10

\skewchar\xptmit='177 \skewchar\xptbmit='177 \skewchar\xptsy='60
\fontdimen16 \xptsy=\the\fontdimen17 \xptsy

\def\xpt{\ifmmode\err@badsizechange\else
     \@mathfontinit
     \textfont0=\xptrm  \scriptfont0=\viiptrm  \scriptscriptfont0=\vptrm
     \textfont1=\xptmit \scriptfont1=\viiptmit \scriptscriptfont1=\vptmit
     \textfont2=\xptsy  \scriptfont2=\viiptsy  \scriptscriptfont2=\vptsy
     \textfont3=\xptex  \scriptfont3=\xptex    \scriptscriptfont3=\xptex
     \textfont\itfam=\xptit
     \scriptfont\itfam=\viiptit
     \scriptscriptfont\itfam=\viiptit
     \textfont\bffam=\xptbf
     \scriptfont\bffam=\viiptbf
     \scriptscriptfont\bffam=\vptbf
     \textfont\bfsfam=\xptbfs
     \scriptfont\bfsfam=\viiptbf
     \scriptscriptfont\bfsfam=\vptbf
     \textfont\bmitfam=\xptbmit
     \scriptfont\bmitfam=\viiptmit
     \scriptscriptfont\bmitfam=\vptmit
     \@fontstyleinit
     \def\rm{\xptrm\fam=\z@}%
     \def\it{\xptit\fam=\itfam}%
     \def\sl{\xptsl}%
     \def\bf{\xptbf\fam=\bffam}%
     \def\tt{\xpttt}%
     \def\ss{\xptss}%
     \def\sc{\xptsc}%
     \def\bfs{\xptbfs\fam=\bfsfam}%
     \def\bmit{\fam=\bmitfam}%
     \def\oldstyle{\xptmit\fam=\@ne}%
     \rm\fi}


\def\getxipt{%
     \font\xiptrm=cmr10  scaled\magstephalf
     \font\xiptmit=cmmi10 scaled\magstephalf
     \font\xiptsy=cmsy10 scaled\magstephalf
     \font\xiptex=cmex10 scaled\magstephalf
     \font\xiptit=cmti10 scaled\magstephalf
     \font\xiptsl=cmsl10 scaled\magstephalf
     \font\xiptbf=cmbx10 scaled\magstephalf
     \font\xipttt=cmtt10 scaled\magstephalf
     \font\xiptss=cmss10 scaled\magstephalf
     \skewchar\xiptmit='177 \skewchar\xiptsy='60
     \fontdimen16 \xiptsy=\the\fontdimen17 \xiptsy}

\def\xipt{\ifmmode\err@badsizechange\else
     \@mathfontinit
     \textfont0=\xiptrm  \scriptfont0=\viiiptrm  \scriptscriptfont0=\viptrm
     \textfont1=\xiptmit \scriptfont1=\viiiptmit \scriptscriptfont1=\viptmit
     \textfont2=\xiptsy  \scriptfont2=\viiiptsy  \scriptscriptfont2=\viptsy
     \textfont3=\xiptex  \scriptfont3=\xptex     \scriptscriptfont3=\xptex
     \textfont\itfam=\xiptit
     \scriptfont\itfam=\viiiptit
     \scriptscriptfont\itfam=\viiptit
     \textfont\bffam=\xiptbf
     \scriptfont\bffam=\viiiptbf
     \scriptscriptfont\bffam=\viptbf
     \@fontstyleinit
     \def\rm{\xiptrm\fam=\z@}%
     \def\it{\xiptit\fam=\itfam}%
     \def\sl{\xiptsl}%
     \def\bf{\xiptbf\fam=\bffam}%
     \def\tt{\xipttt}%
     \def\ss{\xiptss}%
     \def\oldstyle{\xiptmit\fam=\@ne}%
     \rm\fi}


\font\xiiptrm=cmr12 \font\xiiptmit=cmmi12 \font\xiiptsy=cmsy10
scaled\magstep1 \font\xiiptex=cmex10  scaled\magstep1
\font\xiiptit=cmti12 \font\xiiptsl=cmsl12 \font\xiiptbf=cmbx12
\font\xiiptss=cmss12 \font\xiiptsc=cmcsc10 scaled\magstep1
\font\xiiptbfs=cmb10  scaled\magstep1 \font\xiiptbmit=cmmib10
scaled\magstep1

\skewchar\xiiptmit='177 \skewchar\xiiptbmit='177 \skewchar\xiiptsy='60
\fontdimen16 \xiiptsy=\the\fontdimen17 \xiiptsy

\def\xiipt{\ifmmode\err@badsizechange\else
     \@mathfontinit
     \textfont0=\xiiptrm  \scriptfont0=\viiiptrm  \scriptscriptfont0=\viptrm
     \textfont1=\xiiptmit \scriptfont1=\viiiptmit \scriptscriptfont1=\viptmit
     \textfont2=\xiiptsy  \scriptfont2=\viiiptsy  \scriptscriptfont2=\viptsy
     \textfont3=\xiiptex  \scriptfont3=\xptex     \scriptscriptfont3=\xptex
     \textfont\itfam=\xiiptit
     \scriptfont\itfam=\viiiptit
     \scriptscriptfont\itfam=\viiptit
     \textfont\bffam=\xiiptbf
     \scriptfont\bffam=\viiiptbf
     \scriptscriptfont\bffam=\viptbf
     \textfont\bfsfam=\xiiptbfs
     \scriptfont\bfsfam=\viiiptbf
     \scriptscriptfont\bfsfam=\viptbf
     \textfont\bmitfam=\xiiptbmit
     \scriptfont\bmitfam=\viiiptmit
     \scriptscriptfont\bmitfam=\viptmit
     \@fontstyleinit
     \def\rm{\xiiptrm\fam=\z@}%
     \def\it{\xiiptit\fam=\itfam}%
     \def\sl{\xiiptsl}%
     \def\bf{\xiiptbf\fam=\bffam}%
     \def\tt{\xiipttt}%
     \def\ss{\xiiptss}%
     \def\sc{\xiiptsc}%
     \def\bfs{\xiiptbfs\fam=\bfsfam}%
     \def\bmit{\fam=\bmitfam}%
     \def\oldstyle{\xiiptmit\fam=\@ne}%
     \rm\fi}


\def\getxiiipt{%
     \font\xiiiptrm=cmr12  scaled\magstephalf
     \font\xiiiptmit=cmmi12 scaled\magstephalf
     \font\xiiiptsy=cmsy9  scaled\magstep2
     \font\xiiiptit=cmti12 scaled\magstephalf
     \font\xiiiptsl=cmsl12 scaled\magstephalf
     \font\xiiiptbf=cmbx12 scaled\magstephalf
     \font\xiiipttt=cmtt12 scaled\magstephalf
     \font\xiiiptss=cmss12 scaled\magstephalf
     \skewchar\xiiiptmit='177 \skewchar\xiiiptsy='60
     \fontdimen16 \xiiiptsy=\the\fontdimen17 \xiiiptsy}

\def\xiiipt{\ifmmode\err@badsizechange\else
     \@mathfontinit
     \textfont0=\xiiiptrm  \scriptfont0=\xptrm  \scriptscriptfont0=\viiptrm
     \textfont1=\xiiiptmit \scriptfont1=\xptmit \scriptscriptfont1=\viiptmit
     \textfont2=\xiiiptsy  \scriptfont2=\xptsy  \scriptscriptfont2=\viiptsy
     \textfont3=\xivptex   \scriptfont3=\xptex  \scriptscriptfont3=\xptex
     \textfont\itfam=\xiiiptit
     \scriptfont\itfam=\xptit
     \scriptscriptfont\itfam=\viiptit
     \textfont\bffam=\xiiiptbf
     \scriptfont\bffam=\xptbf
     \scriptscriptfont\bffam=\viiptbf
     \@fontstyleinit
     \def\rm{\xiiiptrm\fam=\z@}%
     \def\it{\xiiiptit\fam=\itfam}%
     \def\sl{\xiiiptsl}%
     \def\bf{\xiiiptbf\fam=\bffam}%
     \def\tt{\xiiipttt}%
     \def\ss{\xiiiptss}%
     \def\oldstyle{\xiiiptmit\fam=\@ne}%
     \rm\fi}


\font\xivptrm=cmr12   scaled\magstep1 \font\xivptmit=cmmi12
scaled\magstep1 \font\xivptsy=cmsy10  scaled\magstep2
\font\xivptex=cmex10  scaled\magstep2 \font\xivptit=cmti12
scaled\magstep1 \font\xivptsl=cmsl12  scaled\magstep1
\font\xivptbf=cmbx12  scaled\magstep1
\font\xivptss=cmss12  scaled\magstep1 \font\xivptsc=cmcsc10
scaled\magstep2 \font\xivptbfs=cmb10  scaled\magstep2
\font\xivptbmit=cmmib10 scaled\magstep2

\skewchar\xivptmit='177 \skewchar\xivptbmit='177 \skewchar\xivptsy='60
\fontdimen16 \xivptsy=\the\fontdimen17 \xivptsy

\def\xivpt{\ifmmode\err@badsizechange\else
     \@mathfontinit
     \textfont0=\xivptrm  \scriptfont0=\xptrm  \scriptscriptfont0=\viiptrm
     \textfont1=\xivptmit \scriptfont1=\xptmit \scriptscriptfont1=\viiptmit
     \textfont2=\xivptsy  \scriptfont2=\xptsy  \scriptscriptfont2=\viiptsy
     \textfont3=\xivptex  \scriptfont3=\xptex  \scriptscriptfont3=\xptex
     \textfont\itfam=\xivptit
     \scriptfont\itfam=\xptit
     \scriptscriptfont\itfam=\viiptit
     \textfont\bffam=\xivptbf
     \scriptfont\bffam=\xptbf
     \scriptscriptfont\bffam=\viiptbf
     \textfont\bfsfam=\xivptbfs
     \scriptfont\bfsfam=\xptbfs
     \scriptscriptfont\bfsfam=\viiptbf
     \textfont\bmitfam=\xivptbmit
     \scriptfont\bmitfam=\xptbmit
     \scriptscriptfont\bmitfam=\viiptmit
     \@fontstyleinit
     \def\rm{\xivptrm\fam=\z@}%
     \def\it{\xivptit\fam=\itfam}%
     \def\sl{\xivptsl}%
     \def\bf{\xivptbf\fam=\bffam}%
     \def\tt{\xivpttt}%
     \def\ss{\xivptss}%
     \def\sc{\xivptsc}%
     \def\bfs{\xivptbfs\fam=\bfsfam}%
     \def\bmit{\fam=\bmitfam}%
     \def\oldstyle{\xivptmit\fam=\@ne}%
     \rm\fi}


\font\xviiptrm=cmr17 \font\xviiptmit=cmmi12 scaled\magstep2
\font\xviiptsy=cmsy10 scaled\magstep3 \font\xviiptex=cmex10
scaled\magstep3 \font\xviiptit=cmti12 scaled\magstep2
\font\xviiptbf=cmbx12 scaled\magstep2 \font\xviiptbfs=cmb10
scaled\magstep3

\skewchar\xviiptmit='177 \skewchar\xviiptsy='60 \fontdimen16
\xviiptsy=\the\fontdimen17 \xviiptsy

\def\xviipt{\ifmmode\err@badsizechange\else
     \@mathfontinit
     \textfont0=\xviiptrm  \scriptfont0=\xiiptrm  \scriptscriptfont0=\viiiptrm
     \textfont1=\xviiptmit \scriptfont1=\xiiptmit \scriptscriptfont1=\viiiptmit
     \textfont2=\xviiptsy  \scriptfont2=\xiiptsy  \scriptscriptfont2=\viiiptsy
     \textfont3=\xviiptex  \scriptfont3=\xiiptex  \scriptscriptfont3=\xptex
     \textfont\itfam=\xviiptit
     \scriptfont\itfam=\xiiptit
     \scriptscriptfont\itfam=\viiiptit
     \textfont\bffam=\xviiptbf
     \scriptfont\bffam=\xiiptbf
     \scriptscriptfont\bffam=\viiiptbf
     \textfont\bfsfam=\xviiptbfs
     \scriptfont\bfsfam=\xiiptbfs
     \scriptscriptfont\bfsfam=\viiiptbf
     \@fontstyleinit
     \def\rm{\xviiptrm\fam=\z@}%
     \def\it{\xviiptit\fam=\itfam}%
     \def\bf{\xviiptbf\fam=\bffam}%
     \def\bfs{\xviiptbfs\fam=\bfsfam}%
     \def\oldstyle{\xviiptmit\fam=\@ne}%
     \rm\fi}


\font\xxiptrm=cmr17  scaled\magstep1


\def\xxipt{\ifmmode\err@badsizechange\else
     \@mathfontinit
     \@fontstyleinit
     \def\rm{\xxiptrm\fam=\z@}%
     \rm\fi}


\font\xxvptrm=cmr17  scaled\magstep2


\def\xxvpt{\ifmmode\err@badsizechange\else
     \@mathfontinit
     \@fontstyleinit
     \def\rm{\xxvptrm\fam=\z@}%
     \rm\fi}




\message{Loading jyTeX macros...}

\message{modifications to plain.tex,}


\def\newcount{\alloc@0\count\countdef\insc@unt}
\def\newdimen{\alloc@1\dimen\dimendef\insc@unt}
\def\newskip{\alloc@2\skip\skipdef\insc@unt}
\def\newmuskip{\alloc@3\muskip\muskipdef\@cclvi}
\def\newbox{\alloc@4\box\chardef\insc@unt}
\def\newtoks{\alloc@5\toks\toksdef\@cclvi}
\def\newhelp#1#2{\newtoks#1\global#1\expandafter{\csname#2\endcsname}}
\def\newread{\alloc@6\read\chardef\sixt@@n}
\def\newwrite{\alloc@7\write\chardef\sixt@@n}
\def\newfam{\alloc@8\fam\chardef\sixt@@n}
\def\newinsert#1{\global\advance\insc@unt by\m@ne
     \ch@ck0\insc@unt\count
     \ch@ck1\insc@unt\dimen
     \ch@ck2\insc@unt\skip
     \ch@ck4\insc@unt\box
     \allocationnumber=\insc@unt
     \global\chardef#1=\allocationnumber
     \wlog{\string#1=\string\insert\the\allocationnumber}}
\def\newif#1{\count@\escapechar \escapechar\m@ne
     \expandafter\expandafter\expandafter
          \xdef\@if#1{true}{\let\noexpand#1=\noexpand\iftrue}%
     \expandafter\expandafter\expandafter
          \xdef\@if#1{false}{\let\noexpand#1=\noexpand\iffalse}%
     \global\@if#1{false}\escapechar=\count@}


\newlinechar=`\^^J
\overfullrule=0pt




\let\itfam=\undefined

\let\bffam=\undefined

\count18=3


\chardef\sharps="19


\mathchardef\alpha="710B \mathchardef\beta="710C \mathchardef\gamma="710D
\mathchardef\delta="710E \mathchardef\epsilon="710F
\mathchardef\zeta="7110 \mathchardef\eta="7111 \mathchardef\theta="7112
\mathchardef\iota="7113 \mathchardef\kappa="7114
\mathchardef\lambda="7115 \mathchardef\mu="7116 \mathchardef\nu="7117
\mathchardef\xi="7118 \mathchardef\pi="7119 \mathchardef\rho="711A
\mathchardef\sigma="711B \mathchardef\tau="711C
\mathchardef\upsilon="711D \mathchardef\phi="711E \mathchardef\chi="711F
\mathchardef\psi="7120 \mathchardef\omega="7121
\mathchardef\varepsilon="7122 \mathchardef\vartheta="7123
\mathchardef\varpi="7124 \mathchardef\varrho="7125
\mathchardef\varsigma="7126 \mathchardef\varphi="7127
\mathchardef\imath="717B \mathchardef\jmath="717C \mathchardef\ell="7160
\mathchardef\wp="717D \mathchardef\partial="7140 \mathchardef\flat="715B
\mathchardef\natural="715C \mathchardef\sharp="715D



\def\angle{{\vbox{\ialign{$\m@th\scriptstyle##$\crcr
     \not\mathrel{\mkern14mu}\crcr
     \noalign{\nointerlineskip}
     \mkern2.5mu\leaders\hrule height.34\rp@\hfill\mkern2.5mu\crcr}}}}
\def\vdots{\vbox{\baselineskip4\rp@ \lineskiplimit\z@
     \kern6\rp@\hbox{.}\hbox{.}\hbox{.}}}
\def\ddots{\mathinner{\mkern1mu\raise7\rp@\vbox{\kern7\rp@\hbox{.}}\mkern2mu
     \raise4\rp@\hbox{.}\mkern2mu\raise\rp@\hbox{.}\mkern1mu}}
\def\overrightarrow#1{\vbox{\ialign{##\crcr
     \rightarrowfill\crcr
     \noalign{\kern-\rp@\nointerlineskip}
     $\hfil\displaystyle{#1}\hfil$\crcr}}}
\def\overleftarrow#1{\vbox{\ialign{##\crcr
     \leftarrowfill\crcr
     \noalign{\kern-\rp@\nointerlineskip}
     $\hfil\displaystyle{#1}\hfil$\crcr}}}
\def\overbrace#1{\mathop{\vbox{\ialign{##\crcr
     \noalign{\kern3\rp@}
     \downbracefill\crcr
     \noalign{\kern3\rp@\nointerlineskip}
     $\hfil\displaystyle{#1}\hfil$\crcr}}}\limits}
\def\underbrace#1{\mathop{\vtop{\ialign{##\crcr
     $\hfil\displaystyle{#1}\hfil$\crcr
     \noalign{\kern3\rp@\nointerlineskip}
     \upbracefill\crcr
     \noalign{\kern3\rp@}}}}\limits}
\def\big#1{{\hbox{$\left#1\vbox to8.5\rp@ {}\right.\n@space$}}}
\def\Big#1{{\hbox{$\left#1\vbox to11.5\rp@ {}\right.\n@space$}}}
\def\bigg#1{{\hbox{$\left#1\vbox to14.5\rp@ {}\right.\n@space$}}}
\def\Bigg#1{{\hbox{$\left#1\vbox to17.5\rp@ {}\right.\n@space$}}}
\def\@vereq#1#2{\lower.5\rp@\vbox{\baselineskip\z@skip\lineskip-.5\rp@
     \ialign{$\m@th#1\hfil##\hfil$\crcr#2\crcr=\crcr}}}
\def\rlh@#1{\vcenter{\hbox{\ooalign{\raise2\rp@
     \hbox{$#1\rightharpoonup$}\crcr
     $#1\leftharpoondown$}}}}
\def\bordermatrix#1{\begingroup\m@th
     \setbox\z@\vbox{%
          \def\cr{\crcr\noalign{\kern2\rp@\global\let\cr\endline}}%
          \ialign{$##$\hfil\kern2\rp@\kern\p@renwd
               &\thinspace\hfil$##$\hfil&&\quad\hfil$##$\hfil\crcr
               \omit\strut\hfil\crcr
               \noalign{\kern-\baselineskip}%
               #1\crcr\omit\strut\cr}}%
     \setbox\tw@\vbox{\unvcopy\z@\global\setbox\@ne\lastbox}%
     \setbox\tw@\hbox{\unhbox\@ne\unskip\global\setbox\@ne\lastbox}%
     \setbox\tw@\hbox{$\kern\wd\@ne\kern-\p@renwd\left(\kern-\wd\@ne
          \global\setbox\@ne\vbox{\box\@ne\kern2\rp@}%
          \vcenter{\kern-\ht\@ne\unvbox\z@\kern-\baselineskip}%
          \,\right)$}%
     \null\;\vbox{\kern\ht\@ne\box\tw@}\endgroup}
\def\endinsert{\egroup
     \if@mid\dimen@\ht\z@
          \advance\dimen@\dp\z@
          \advance\dimen@12\rp@
          \advance\dimen@\pagetotal
          \ifdim\dimen@>\pagegoal\@midfalse\p@gefalse\fi
     \fi
     \if@mid\bigskip\box\z@
          \bigbreak
     \else\insert\topins{\penalty100 \splittopskip\z@skip
               \splitmaxdepth\maxdimen\floatingpenalty\z@
               \ifp@ge\dimen@\dp\z@
                    \vbox to\vsize{\unvbox\z@\kern-\dimen@}%
               \else\box\z@\nobreak\bigskip
               \fi}%
     \fi
     \endgroup}


\def\cases#1{\left\{\,\vcenter{\m@th
     \ialign{$##\hfil$&\quad##\hfil\crcr#1\crcr}}\right.}
\def\matrix#1{\null\,\vcenter{\m@th
     \ialign{\hfil$##$\hfil&&\quad\hfil$##$\hfil\crcr
          \mathstrut\crcr
          \noalign{\kern-\baselineskip}
          #1\crcr
          \mathstrut\crcr
          \noalign{\kern-\baselineskip}}}\,}


\newif\ifraggedbottom

\def\raggedbottom{\ifraggedbottom\else
     \advance\topskip by\z@ plus60pt \raggedbottomtrue\fi}%
\def\normalbottom{\ifraggedbottom
     \advance\topskip by\z@ plus-60pt \raggedbottomfalse\fi}

\message{hacks,}


\toksdef\toks@i=1 \toksdef\toks@ii=2


\def\TeX{T\kern-.1667em \lower.5ex \hbox{E}\kern-.125em X\null}
\def\jyTeX{{\leavevmode
     \raise.587ex \hbox{\it\j}\kern-.1em \lower.048ex \hbox{\it y}\kern-.12em
     \TeX}}

\let\then=\iftrue
\def\ifnoarg#1\then{\def\hack@{#1}\ifx\hack@\empty}
\def\ifundefined#1\then{%
     \expandafter\ifx\csname\expandafter\blank\string#1\endcsname\relax}
\def\useif#1\then{\csname#1\endcsname}
\def\usename#1{\csname#1\endcsname}
\def\useafter#1#2{\expandafter#1\csname#2\endcsname}

\long\def\loop#1\repeat{\def\@iterate{#1\expandafter\@iterate\fi}\@iterate
     \let\@iterate=\relax}

\let\TeXend=\end
\def\begin#1{\begingroup\def\@@blockname{#1}\usename{begin#1}}
\def\end#1{\usename{end#1}\def\hack@{#1}%
     \ifx\@@blockname\hack@
          \endgroup
     \else\err@badgroup\hack@\@@blockname
     \fi}
\def\@@blockname{}

\def\defaultoption[#1]#2{%
     \def\hack@{\ifx\hack@ii[\toks@={#2}\else\toks@={#2[#1]}\fi\the\toks@}%
     \futurelet\hack@ii\hack@}

\def\markup#1{\let\@@marksf=\empty
     \ifhmode\edef\@@marksf{\spacefactor=\the\spacefactor\relax}\/\fi
     ${}^{\hbox{\subscriptfonts#1}}$\@@marksf}


\newtoks\shortyear
\newtoks\militaryhour
\newtoks\standardhour
\newtoks\minute
\newtoks\amorpm

\def\settime{\count@=\time\divide\count@ by60
     \militaryhour=\expandafter{\number\count@}%
     {\multiply\count@ by-60 \advance\count@ by\time
          \xdef\hack@{\ifnum\count@<10 0\fi\number\count@}}%
     \minute=\expandafter{\hack@}%
     \ifnum\count@<12
          \amorpm={am}
     \else\amorpm={pm}
          \ifnum\count@>12 \advance\count@ by-12 \fi
     \fi
     \standardhour=\expandafter{\number\count@}%
     \def\hack@19##1##2{\shortyear={##1##2}}%
          \expandafter\hack@\the\year}

\def\monthword#1{%
     \ifcase#1
          $\bullet$\err@badcountervalue{monthword}%
          \or January\or February\or March\or April\or May\or June%
          \or July\or August\or September\or October\or November\or December%
     \else$\bullet$\err@badcountervalue{monthword}%
     \fi}

\def\monthabbr#1{%
     \ifcase#1
          $\bullet$\err@badcountervalue{monthabbr}%
          \or Jan\or Feb\or Mar\or Apr\or May\or Jun%
          \or Jul\or Aug\or Sep\or Oct\or Nov\or Dec%
     \else$\bullet$\err@badcountervalue{monthabbr}%
     \fi}

\def\militarytime{\the\militaryhour:\the\minute}
\def\standardtime{\the\standardhour:\the\minute}


\def\@setnumstyle#1#2{\expandafter\global\expandafter\expandafter
     \expandafter\let\expandafter\expandafter
     \csname @\expandafter\blank\string#1style\endcsname
     \csname#2\endcsname}
\def\numstyle#1{\usename{@\expandafter\blank\string#1style}#1}
\def\ifblank#1\then{\useafter\ifx{@\expandafter\blank\string#1}\blank}

\def\blank#1{}

\def\Roman#1{\expandafter\uppercase\expandafter{\romannumeral#1}}
\def\alphabetic#1{%
     \ifcase#1
          $\bullet$\err@badcountervalue{alphabetic}%
          \or a\or b\or c\or d\or e\or f\or g\or h\or i\or j\or k\or l\or m%
          \or n\or o\or p\or q\or r\or s\or t\or u\or v\or w\or x\or y\or z%
     \else$\bullet$\err@badcountervalue{alphabetic}%
     \fi}
\def\Alphabetic#1{\expandafter\uppercase\expandafter{\alphabetic{#1}}}
\def\symbols#1{%
     \ifcase#1
          $\bullet$\err@badcountervalue{symbols}%
          \or*\or\dag\or\ddag\or\S\or$\|$%
          \or**\or\dag\dag\or\ddag\ddag\or\S\S\or$\|\|$%
     \else$\bullet$\err@badcountervalue{symbols}%
     \fi}


\catcode`\^^?=13 \def^^?{\relax}

\def\trimleading#1\to#2{\edef#2{#1}%
     \expandafter\@trimleading\expandafter#2#2^^?^^?}
\def\@trimleading#1#2#3^^?{\ifx#2^^?\def#1{}\else\def#1{#2#3}\fi}

\def\trimtrailing#1\to#2{\edef#2{#1}%
     \expandafter\@trimtrailing\expandafter#2#2^^? ^^?\relax}
\def\@trimtrailing#1#2 ^^?#3{\ifx#3\relax\toks@={}%
     \else\def#1{#2}\toks@={\trimtrailing#1\to#1}\fi
     \the\toks@}

\def\trim#1\to#2{\trimleading#1\to#2\trimtrailing#2\to#2}

\catcode`\^^?=15


\long\def\additemL#1\to#2{\toks@={\^^\{#1}}\toks@ii=\expandafter{#2}%
     \xdef#2{\the\toks@\the\toks@ii}}

\long\def\additemR#1\to#2{\toks@={\^^\{#1}}\toks@ii=\expandafter{#2}%
     \xdef#2{\the\toks@ii\the\toks@}}

\def\getitemL#1\to#2{\expandafter\@getitemL#1\hack@#1#2}
\def\@getitemL\^^\#1#2\hack@#3#4{\def#4{#1}\def#3{#2}}

\message{font macros,}


\newdimen\rp@
\newcount\@@sizeindex \@@sizeindex=0
\newcount\@@factori
\newcount\@@factorii
\newcount\@@factoriii
\newcount\@@factoriv

\countdef\maxfam=18
\newfam\itfam
\newfam\bffam
\newfam\bfsfam
\newfam\bmitfam

\def\@mathfontinit{\count@=4
     \loop\textfont\count@=\nullfont
          \scriptfont\count@=\nullfont
          \scriptscriptfont\count@=\nullfont
          \ifnum\count@<\maxfam\advance\count@ by\@ne
     \repeat}

\def\@fontstyleinit{%
     \def\it{\err@fontnotavailable\it}%
     \def\bf{\err@fontnotavailable\bf}%
     \def\bfs{\err@bfstobf}%
     \def\bmit{\err@fontnotavailable\bmit}%
     \def\sc{\err@fontnotavailable\sc}%
     \def\sl{\err@sltoit}%
     \def\ss{\err@fontnotavailable\ss}%
     \def\tt{\err@fontnotavailable\tt}}

\def\@parameterinit#1{\rm\rp@=.1em \@getscaling{#1}%
     \let\^^\=\@doscaling\scalingskipslist
     \setbox\strutbox=\hbox{\vrule
          height.708\baselineskip depth.292\baselineskip width\z@}}

\def\@getfactor#1#2#3#4{\@@factori=#1 \@@factorii=#2
     \@@factoriii=#3 \@@factoriv=#4}

\def\@getscaling#1{\count@=#1 \advance\count@ by-\@@sizeindex\@@sizeindex=#1
     \ifnum\count@<0
          \let\@mulordiv=\divide
          \let\@divormul=\multiply
          \multiply\count@ by\m@ne
     \else\let\@mulordiv=\multiply
          \let\@divormul=\divide
     \fi
     \edef\@@scratcha{\ifcase\count@                {1}{1}{1}{1}\or
          {1}{7}{23}{3}\or     {2}{5}{3}{1}\or      {9}{89}{13}{1}\or
          {6}{25}{6}{1}\or     {8}{71}{14}{1}\or    {6}{25}{36}{5}\or
          {1}{7}{53}{4}\or     {12}{125}{108}{5}\or {3}{14}{53}{5}\or
          {6}{41}{17}{1}\or    {13}{31}{13}{2}\or   {9}{107}{71}{2}\or
          {11}{139}{124}{3}\or {1}{6}{43}{2}\or     {10}{107}{42}{1}\or
          {1}{5}{43}{2}\or     {5}{69}{65}{1}\or    {11}{97}{91}{2}\fi}%
     \expandafter\@getfactor\@@scratcha}

\def\@doscaling#1{\@mulordiv#1by\@@factori\@divormul#1by\@@factorii
     \@mulordiv#1by\@@factoriii\@divormul#1by\@@factoriv}


\newskip\headskip
\newskip\footskip

\def\typesize=#1pt{\count@=#1 \advance\count@ by-10
     \ifcase\count@
          \@setsizex\or\err@badtypesize\or
          \@setsizexii\or\err@badtypesize\or
          \@setsizexiv
     \else\err@badtypesize
     \fi}

\def\@setsizex{\getixpt
     \def\subsubscriptfonts{\vpt}%
          \def\subsubscriptsize{\vpt\@parameterinit{-8}}%
     \def\subscriptfonts{\viipt}\def\subscriptsize{\viipt\@parameterinit{-4}}%
     \def\footnotefonts{\viiipt}\def\footnotesize{\viiipt\@parameterinit{-2}}%
     \def\smallfonts{\ixpt}\def\smallsize{\ixpt\@parameterinit{-1}}%
     \def\normalfonts{\xpt}\def\normalsize{\xpt\@parameterinit{0}}%
     \def\bigfonts{\xiipt}\def\bigsize{\xiipt\@parameterinit{2}}%
     \def\Bigfonts{\xivpt}\def\Bigsize{\xivpt\@parameterinit{4}}%
     \def\biggfonts{\xviipt}\def\biggsize{\xviipt\@parameterinit{6}}%
     \def\Biggfonts{\xxipt}\def\Biggsize{\xxipt\@parameterinit{8}}%
     \def\tinyfonts{\vpt}\def\tinysize{\vpt\@parameterinit{-8}}%
     \def\HUGEFONTS{\xxvpt}\def\HUGESIZE{\xxvpt\@parameterinit{10}}%
     \normalsize\fixedskipslist}

\def\@setsizexii{\getxipt
     \def\subsubscriptfonts{\vipt}%
          \def\subsubscriptsize{\vipt\@parameterinit{-6}}%
     \def\subscriptfonts{\viiipt}%
          \def\subscriptsize{\viiipt\@parameterinit{-2}}%
     \def\footnotefonts{\xpt}\def\footnotesize{\xpt\@parameterinit{0}}%
     \def\smallfonts{\xipt}\def\smallsize{\xipt\@parameterinit{1}}%
     \def\normalfonts{\xiipt}\def\normalsize{\xiipt\@parameterinit{2}}%
     \def\bigfonts{\xivpt}\def\bigsize{\xivpt\@parameterinit{4}}%
     \def\Bigfonts{\xviipt}\def\Bigsize{\xviipt\@parameterinit{6}}%
     \def\biggfonts{\xxipt}\def\biggsize{\xxipt\@parameterinit{8}}%
     \def\Biggfonts{\xxvpt}\def\Biggsize{\xxvpt\@parameterinit{10}}%
     \def\tinyfonts{\vpt}\def\tinysize{\vpt\@parameterinit{-8}}%
     \def\HUGEFONTS{\xxvpt}\def\HUGESIZE{\xxvpt\@parameterinit{10}}%
     \normalsize\fixedskipslist}

\def\@setsizexiv{\getxiiipt
     \def\subsubscriptfonts{\viipt}%
          \def\subsubscriptsize{\viipt\@parameterinit{-4}}%
     \def\subscriptfonts{\xpt}\def\subscriptsize{\xpt\@parameterinit{0}}%
     \def\footnotefonts{\xiipt}\def\footnotesize{\xiipt\@parameterinit{2}}%
     \def\smallfonts{\xiiipt}\def\smallsize{\xiiipt\@parameterinit{3}}%
     \def\normalfonts{\xivpt}\def\normalsize{\xivpt\@parameterinit{4}}%
     \def\bigfonts{\xviipt}\def\bigsize{\xviipt\@parameterinit{6}}%
     \def\Bigfonts{\xxipt}\def\Bigsize{\xxipt\@parameterinit{8}}%
     \def\biggfonts{\xxvpt}\def\biggsize{\xxvpt\@parameterinit{10}}%
     \def\Biggfonts{\err@sizetoolarge\Biggfonts\HUGEFONTS}%
          \def\Biggsize{\err@sizetoolarge\Biggsize\HUGESIZE}%
     \def\tinyfonts{\vpt}\def\tinysize{\vpt\@parameterinit{-8}}%
     \def\HUGEFONTS{\xxvpt}\def\HUGESIZE{\xxvpt\@parameterinit{10}}%
     \normalsize\fixedskipslist}

\def\subsubscriptfonts{\vpt} \def\subsubscriptsize{\vpt\@parameterinit{-8}}
\def\subscriptfonts{\viipt}  \def\subscriptsize{\viipt\@parameterinit{-4}}
\def\footnotefonts{\viiipt}  \def\footnotesize{\viiipt\@parameterinit{-2}}
\def\smallfonts{\err@sizenotavailable\smallfonts}
                             \def\smallsize{\ixpt\@parameterinit{-1}}
\def\normalfonts{\xpt}       \def\normalsize{\xpt\@parameterinit{0}}
\def\bigfonts{\xiipt}        \def\bigsize{\xiipt\@parameterinit{2}}
\def\Bigfonts{\xivpt}        \def\Bigsize{\xivpt\@parameterinit{4}}
\def\biggfonts{\xviipt}      \def\biggsize{\xviipt\@parameterinit{6}}
\def\Biggfonts{\xxipt}       \def\Biggsize{\xxipt\@parameterinit{8}}
\def\tinyfonts{\vpt}         \def\tinysize{\vpt\@parameterinit{-8}}
\def\HUGEFONTS{\xxvpt}       \def\HUGESIZE{\xxvpt\@parameterinit{10}}

\message{document layout,}


\newtoks\everyoutput \everyoutput={}
\newdimen\depthofpage
\newcount\pagenum \pagenum=0

\newdimen\oddtopmargin  \newdimen\eventopmargin
\newdimen\oddleftmargin \newdimen\evenleftmargin
\newtoks\oddhead        \newtoks\evenhead
\newtoks\oddfoot        \newtoks\evenfoot

\def\topmargin{\afterassignment\@seteventop\oddtopmargin}
\def\leftmargin{\afterassignment\@setevenleft\oddleftmargin}
\def\head{\afterassignment\@setevenhead\oddhead}
\def\foot{\afterassignment\@setevenfoot\oddfoot}

\def\@seteventop{\eventopmargin=\oddtopmargin}
\def\@setevenleft{\evenleftmargin=\oddleftmargin}
\def\@setevenhead{\evenhead=\oddhead}
\def\@setevenfoot{\evenfoot=\oddfoot}

\def\pagenumstyle#1{\@setnumstyle\pagenum{#1}}

\newif\ifdraft
\def\draft{\drafttrue\leftmargin=.5in \overfullrule=5pt }

\def\outputstyle#1{\global\expandafter\let\expandafter
          \@outputstyle\csname#1output\endcsname
     \usename{#1setup}}

\output={\@outputstyle}

\def\normaloutput{\the\everyoutput
     \global\advance\pagenum by\@ne
     \ifodd\pagenum
          \voffset=\oddtopmargin \hoffset=\oddleftmargin
     \else\voffset=\eventopmargin \hoffset=\evenleftmargin
     \fi
     \advance\voffset by-1in  \advance\hoffset by-1in
     \count0=\pagenum
     \expandafter\shipout\pagebox
     \ifnum\outputpenalty>-\@MM\else\dosupereject\fi}

\newdimen\fullhsize
\newbox\leftpage
\newcount\leftpagenum
\newcount\outputpagenum \outputpagenum=0
\let\leftorright=L

\def\twoupoutput{\the\everyoutput
     \global\advance\pagenum by\@ne
     \if L\leftorright
          \global\setbox\leftpage=\leftline{\pagebox}%
          \global\leftpagenum=\pagenum
          \global\let\leftorright=R%
     \else\global\advance\outputpagenum by\@ne
          \ifodd\outputpagenum
               \voffset=\oddtopmargin \hoffset=\oddleftmargin
          \else\voffset=\eventopmargin \hoffset=\evenleftmargin
          \fi
          \advance\voffset by-1in  \advance\hoffset by-1in
          \count0=\leftpagenum \count1=\pagenum
          \shipout\vbox{\hbox to\fullhsize
               {\box\leftpage\hfil\leftline{\pagebox}}}%
          \global\let\leftorright=L%
     \fi
     \ifnum\outputpenalty>-\@MM
     \else\dosupereject
          \if R\leftorright
               \globaldefs=\@ne\head={\hfil}\foot={\hfil}\globaldefs=\z@
               \null\newpage
          \fi
     \fi}

\def\pagebox{\vbox{\makeheadline\pagebody\makefootline}}

\def\makeheadline{%
     \vbox to\z@{\baselinestretch=\@m
          \vskip\topskip\vskip-.708\baselineskip\vskip-\headskip
          \line{\vbox to\ht\strutbox{}%
               \ifodd\pagenum\the\oddhead\else\the\evenhead\fi}%
          \vss}%
     \nointerlineskip}

\def\pagebody{\vbox to\vsize{%
     \boxmaxdepth\maxdepth
     \ifvoid\topins\else\unvbox\topins\fi
     \depthofpage=\dp255
     \unvbox255
     \ifraggedbottom\kern-\depthofpage\vfil\fi
     \ifvoid\footins
     \else\vskip\skip\footins
          \footnoterule
          \unvbox\footins
          \vskip-\footnoteskip
     \fi}}

\def\makefootline{\baselineskip=\footskip
     \line{\ifodd\pagenum\the\oddfoot\else\the\evenfoot\fi}}


\newskip\abovechapterskip
\newskip\belowchapterskip
\newskip\abovesectionskip
\newskip\belowsectionskip
\newskip\abovesubsectionskip
\newskip\belowsubsectionskip

\def\chapterstyle#1{\global\expandafter\let\expandafter\@chapterstyle
     \csname#1text\endcsname}
\def\sectionstyle#1{\global\expandafter\let\expandafter\@sectionstyle
     \csname#1text\endcsname}
\def\subsectionstyle#1{\global\expandafter\let\expandafter\@subsectionstyle
     \csname#1text\endcsname}

\def\chapter#1{%
     \ifdim\lastskip=17sp \else\chapterbreak\vskip\abovechapterskip\fi
     \@chapterstyle{\ifblank\chapternumstyle\then
          \else\newchapternum=\next\chapternumformat\ \fi#1}%
     \nobreak\vskip\belowchapterskip\vskip17sp }

\def\section#1{%
     \ifdim\lastskip=17sp \else\sectionbreak\vskip\abovesectionskip\fi
     \@sectionstyle{\ifblank\sectionnumstyle\then
          \else\newsectionnum=\next\sectionnumformat\ \fi#1}%
     \nobreak\vskip\belowsectionskip\vskip17sp }

\def\subsection#1{%
     \ifdim\lastskip=17sp \else\subsectionbreak\vskip\abovesubsectionskip\fi
     \@subsectionstyle{\ifblank\subsectionnumstyle\then
          \else\newsubsectionnum=\next\subsectionnumformat\ \fi#1}%
     \nobreak\vskip\belowsubsectionskip\vskip17sp }


\let\TeXunderline=\underline
\let\TeXoverline=\overline
\def\underline#1{\relax\ifmmode\TeXunderline{#1}\else
     $\TeXunderline{\hbox{#1}}$\fi}
\def\overline#1{\relax\ifmmode\TeXoverline{#1}\else
     $\TeXoverline{\hbox{#1}}$\fi}

\def\baselinestretch{\afterassignment\@baselinestretch\count@}
\def\@baselinestretch{\baselineskip=\normalbaselineskip
     \divide\baselineskip by\@m\baselineskip=\count@\baselineskip
     \setbox\strutbox=\hbox{\vrule
          height.708\baselineskip depth.292\baselineskip width\z@}%
     \bigskipamount=\the\baselineskip
          plus.25\baselineskip minus.25\baselineskip
     \medskipamount=.5\baselineskip
          plus.125\baselineskip minus.125\baselineskip
     \smallskipamount=.25\baselineskip
          plus.0625\baselineskip minus.0625\baselineskip}

\def\\{\ifhmode\ifnum\lastpenalty=-\@M\else\hfil\penalty-\@M\fi\fi
     \ignorespaces}
\def\newpage{\vfil\break}

\def\lefttext#1{\par{\@text\leftskip=\z@\rightskip=\centering
     \noindent#1\par}}
\def\righttext#1{\par{\@text\leftskip=\centering\rightskip=\z@
     \noindent#1\par}}
\def\centertext#1{\par{\@text\leftskip=\centering\rightskip=\centering
     \noindent#1\par}}
\def\@text{\parindent=\z@ \parfillskip=\z@ \everypar={}%
     \spaceskip=.3333em \xspaceskip=.5em
     \def\\{\ifhmode\ifnum\lastpenalty=-\@M\else\penalty-\@M\fi\fi
          \ignorespaces}}

\def\beginleft{\par\@text\leftskip=\z@ \rightskip=\centering}
     
\def\beginright{\par\@text\leftskip=\centering\rightskip=\z@ }
     
\def\begincenter{\par\@text\leftskip=\centering\rightskip=\centering}

\def\beginnarrow{\defaultoption[\parindent]\@beginnarrow}
\def\@beginnarrow[#1]{\par\advance\leftskip by#1\advance\rightskip by#1}

\begingroup
\catcode`\[=1 \catcode`\{=11 \gdef\beginignore[\endgroup\bgroup
     \catcode`\e=0 \catcode`\\=12 \catcode`\{=11 \catcode`\f=12 \let\or=\relax
     \let\nd{ignor=\fi \let\}=\egroup
     \iffalse}
\endgroup

\long\def\marginnote#1{\leavevmode
     \edef\@marginsf{\spacefactor=\the\spacefactor\relax}%
     \ifdraft\strut\vadjust{%
          \hbox to\z@{\hskip\hsize\hskip.1in
               \vbox to\z@{\vskip-\dp\strutbox
                    \marginnoteformat
                    \vskip-\ht\strutbox
                    \noindent\strut#1\par
                    \vss}%
               \hss}}%
     \fi
     \@marginsf}


\newtoks\everybye \everybye={\par\vfil}
\outer\def\bye{\the\everybye
     \footnotecheck
     \prelabelcheck
     \streamcheck
     \supereject
     \TeXend}

\message{footnotes,}

\newcount\footnotenum \footnotenum=0
\newskip\footnoteskip
\let\@footnotelist=\empty

\def\footnotenumstyle#1{\@setnumstyle\footnotenum{#1}%
     \useafter\ifx{@footnotenumstyle}\symbols
          \global\let\@footup=\empty
     \else\global\let\@footup=\markup
     \fi}

\def\footnote{\footnotecheck\defaultoption[]\@footnote}
\def\@footnote[#1]{\@footnotemark[#1]\@footnotetext}

\def\footnotemark{\defaultoption[]\@footnotemark}
\def\@footnotemark[#1]{\let\@footsf=\empty
     \ifhmode\edef\@footsf{\spacefactor=\the\spacefactor\relax}\/\fi
     \ifnoarg#1\then
          \global\advance\footnotenum by\@ne
          \@footup{\footnotenumformat}%
          \edef\@@foota{\footnotenum=\the\footnotenum\relax}%
          \expandafter\additemR\expandafter\@footup\expandafter
               {\@@foota\footnotenumformat}\to\@footnotelist
          \global\let\@footnotelist=\@footnotelist
     \else\markup{#1}%
          \additemR\markup{#1}\to\@footnotelist
          \global\let\@footnotelist=\@footnotelist
     \fi
     \@footsf}

\def\footnotetext{%
     \ifx\@footnotelist\empty\err@extrafootnotetext\else\@footnotetext\fi}
\def\@footnotetext{%
     \getitemL\@footnotelist\to\@@foota
     \global\let\@footnotelist=\@footnotelist
     \insert\footins\bgroup
     \footnoteformat
     \splittopskip=\ht\strutbox\splitmaxdepth=\dp\strutbox
     \interlinepenalty=\interfootnotelinepenalty\floatingpenalty=\@MM
     \noindent\llap{\@@foota}\strut
     \bgroup\aftergroup\@footnoteend
     \let\@@scratcha=}
\def\@footnoteend{\strut\par\vskip\footnoteskip\egroup}

\def\footnoterule{\normalfonts
     \kern-.3em \hrule width2in height.04em \kern .26em }

\def\footnotecheck{%
     \ifx\@footnotelist\empty
     \else\err@extrafootnotemark
          \global\let\@footnotelist=\empty
     \fi}

\message{labels,}

\let\@@labeldef=\xdef
\newif\if@labelfile
\newwrite\@labelfile
\let\@prelabellist=\empty

\def\label#1#2{\trim#1\to\@@labarg\edef\@@labtext{#2}%
     \edef\@@labname{lab@\@@labarg}%
     \useafter\ifundefined\@@labname\then\else\@yeslab\fi
     \useafter\@@labeldef\@@labname{#2}%
     \ifstreaming
          \expandafter\toks@\expandafter\expandafter\expandafter
               {\csname\@@labname\endcsname}%
          \immediate\write\streamout{\noexpand\label{\@@labarg}{\the\toks@}}%
     \fi}
\def\@yeslab{%
     \useafter\ifundefined{if\@@labname}\then
          \err@labelredef\@@labarg
     \else\useif{if\@@labname}\then
               \err@labelredef\@@labarg
          \else\global\usename{\@@labname true}%
               \useafter\ifundefined{pre\@@labname}\then
               \else\useafter\ifx{pre\@@labname}\@@labtext
                    \else\err@badlabelmatch\@@labarg
                    \fi
               \fi
               \if@labelfile
               \else\global\@labelfiletrue
                    \immediate\write\sixt@@n{--> Creating file \jobname.lab}%
                    \immediate\openout\@labelfile=\jobname.lab
               \fi
               \immediate\write\@labelfile
                    {\noexpand\prelabel{\@@labarg}{\@@labtext}}%
          \fi
     \fi}

\def\putlab#1{\trim#1\to\@@labarg\edef\@@labname{lab@\@@labarg}%
     \useafter\ifundefined\@@labname\then\@nolab\else\usename\@@labname\fi}
\def\@nolab{%
     \useafter\ifundefined{pre\@@labname}\then
          \undefinedlabelformat
          \err@needlabel\@@labarg
          \useafter\xdef\@@labname{\undefinedlabelformat}%
     \else\usename{pre\@@labname}%
          \useafter\xdef\@@labname{\usename{pre\@@labname}}%
     \fi
     \useafter\newif{if\@@labname}%
     \expandafter\additemR\@@labarg\to\@prelabellist}

\def\prelabel#1{\useafter\gdef{prelab@#1}}

\def\ifundefinedlabel#1\then{%
     \expandafter\ifx\csname lab@#1\endcsname\relax}
\def\useiflab#1\then{\csname iflab@#1\endcsname}

\def\prelabelcheck{{%
     \def\^^\##1{\useiflab{##1}\then\else\err@undefinedlabel{##1}\fi}%
     \@prelabellist}}

\message{equation numbering,}

\newcount\chapternum
\newcount\sectionnum
\newcount\subsectionnum
\newcount\equationnum
\newcount\subequationnum
\newcount\figurenum
\newcount\subfigurenum
\newcount\tablenum
\newcount\subtablenum

\newif\if@subeqncount
\newif\if@subfigcount
\newif\if@subtblcount

\def\newchapternum{\newsectionnum=\z@\@resetnum\chapternum}
\def\newsectionnum{\newsubsectionnum=\z@\@resetnum\sectionnum}
\def\newsubsectionnum{\newequationnum=\z@\newfigurenum=\z@\newtablenum=\z@
     \@resetnum\subsectionnum}
\def\newequationnum{\newsubequationnum=\z@\@resetnum\equationnum}
\def\newsubequationnum{\@resetnum\subequationnum}
\def\newfigurenum{\newsubfigurenum=\z@\@resetnum\figurenum}
\def\newsubfigurenum{\@resetnum\subfigurenum}
\def\newtablenum{\newsubtablenum=\z@\@resetnum\tablenum}
\def\newsubtablenum{\@resetnum\subtablenum}

\def\@resetnum#1{\global\advance#1by1 \edef\next{\the#1\relax}\global#1}

\newchapternum=0

\def\chapternumstyle#1{\@setnumstyle\chapternum{#1}}
\def\sectionnumstyle#1{\@setnumstyle\sectionnum{#1}}
\def\subsectionnumstyle#1{\@setnumstyle\subsectionnum{#1}}
\def\equationnumstyle#1{\@setnumstyle\equationnum{#1}}
\def\subequationnumstyle#1{\@setnumstyle\subequationnum{#1}%
     \ifblank\subequationnumstyle\then\global\@subeqncountfalse\fi
     \ignorespaces}
\def\figurenumstyle#1{\@setnumstyle\figurenum{#1}}
\def\subfigurenumstyle#1{\@setnumstyle\subfigurenum{#1}%
     \ifblank\subfigurenumstyle\then\global\@subfigcountfalse\fi
     \ignorespaces}
\def\tablenumstyle#1{\@setnumstyle\tablenum{#1}}
\def\subtablenumstyle#1{\@setnumstyle\subtablenum{#1}%
     \ifblank\subtablenumstyle\then\global\@subtblcountfalse\fi
     \ignorespaces}

\def\eqnlabel#1{%
     \if@subeqncount
          \newsubequationnum=\next
     \else\newequationnum=\next
          \ifblank\subequationnumstyle\then
          \else\global\@subeqncounttrue
               \newsubequationnum=\@ne
          \fi
     \fi
     \label{#1}{\puteqnformat}(\puteqn{#1})%
     \ifdraft\rlap{\hskip.1in{\tt#1}}\fi}

\let\puteqn=\putlab

\def\equation#1#2{\useafter\gdef{eqn@#1}{#2\eqno\eqnlabel{#1}}}
\def\Equation#1{\useafter\gdef{eqn@#1}}

\def\putequation#1{\useafter\ifundefined{eqn@#1}\then
     \err@undefinedeqn{#1}\else\usename{eqn@#1}\fi}

\def\eqnseriesstyle#1{\gdef\@eqnseriesstyle{#1}}
\def\begineqnseries{\subequationnumstyle{\@eqnseriesstyle}%
     \defaultoption[]\@begineqnseries}
\def\@begineqnseries[#1]{\edef\@@eqnname{#1}}
\def\endeqnseries{\subequationnumstyle{blank}%
     \expandafter\ifnoarg\@@eqnname\then
     \else\label\@@eqnname{\puteqnformat}%
     \fi
     \aftergroup\ignorespaces}

\def\figlabel#1{%
     \if@subfigcount
          \newsubfigurenum=\next
     \else\newfigurenum=\next
          \ifblank\subfigurenumstyle\then
          \else\global\@subfigcounttrue
               \newsubfigurenum=\@ne
          \fi
     \fi
     \label{#1}{\putfigformat}\putfig{#1}%
     {\def\marginnoteformat{\tt}\marginnote{#1}}}

\let\putfig=\putlab

\def\figseriesstyle#1{\gdef\@figseriesstyle{#1}}
\def\beginfigseries{\subfigurenumstyle{\@figseriesstyle}%
     \defaultoption[]\@beginfigseries}
\def\@beginfigseries[#1]{\edef\@@figname{#1}}
\def\endfigseries{\subfigurenumstyle{blank}%
     \expandafter\ifnoarg\@@figname\then
     \else\label\@@figname{\putfigformat}%
     \fi
     \aftergroup\ignorespaces}

\def\tbllabel#1{%
     \if@subtblcount
          \newsubtablenum=\next
     \else\newtablenum=\next
          \ifblank\subtablenumstyle\then
          \else\global\@subtblcounttrue
               \newsubtablenum=\@ne
          \fi
     \fi
     \label{#1}{\puttblformat}\puttbl{#1}%
     {\def\marginnoteformat{\tt}\marginnote{#1}}}

\let\puttbl=\putlab

\def\tblseriesstyle#1{\gdef\@tblseriesstyle{#1}}
\def\begintblseries{\subtablenumstyle{\@tblseriesstyle}%
     \defaultoption[]\@begintblseries}
\def\@begintblseries[#1]{\edef\@@tblname{#1}}
\def\endtblseries{\subtablenumstyle{blank}%
     \expandafter\ifnoarg\@@tblname\then
     \else\label\@@tblname{\puttblformat}%
     \fi
     \aftergroup\ignorespaces}

\message{reference numbering,}

\newcount\referencenum \referencenum=0
\newcount\@@prerefcount \@@prerefcount=0
\newcount\@@thisref
\newcount\@@lastref
\newcount\@@loopref
\newcount\@@refseq
\newdimen\refnumindent
\let\@undefreflist=\empty

\def\referencenumstyle#1{\@setnumstyle\referencenum{#1}}

\def\referencestyle#1{\usename{@ref#1}}

\def\@refsequential{%
     \gdef\@refpredef##1{\global\advance\referencenum by\@ne
          \let\^^\=0\label{##1}{\^^\{\the\referencenum}}%
          \useafter\gdef{ref@\the\referencenum}{{##1}{\undefinedlabelformat}}}%
     \gdef\@reference##1##2{%
          \ifundefinedlabel##1\then
          \else\def\^^\####1{\global\@@thisref=####1\relax}\putlab{##1}%
               \useafter\gdef{ref@\the\@@thisref}{{##1}{##2}}%
          \fi}%
     \gdef\endputreferences{%
          \loop\ifnum\@@loopref<\referencenum
                    \advance\@@loopref by\@ne
                    \expandafter\expandafter\expandafter\@printreference
                         \csname ref@\the\@@loopref\endcsname
          \repeat
          \par}}

\def\@refpreordered{%
     \gdef\@refpredef##1{\global\advance\referencenum by\@ne
          \additemR##1\to\@undefreflist}%
     \gdef\@reference##1##2{%
          \ifundefinedlabel##1\then
          \else\global\advance\@@loopref by\@ne
               {\let\^^\=0\label{##1}{\^^\{\the\@@loopref}}}%
               \@printreference{##1}{##2}%
          \fi}
     \gdef\endputreferences{%
          \def\^^\####1{\useiflab{####1}\then
               \else\reference{####1}{\undefinedlabelformat}\fi}%
          \@undefreflist
          \par}}

\def\beginprereferences{\par
     \def\reference##1##2{\global\advance\referencenum by1\@ne
          \let\^^\=0\label{##1}{\^^\{\the\referencenum}}%
          \useafter\gdef{ref@\the\referencenum}{{##1}{##2}}}}
\def\endprereferences{\global\@@prerefcount=\the\referencenum\par}

\def\beginputreferences{\par
     \refnumindent=\z@\@@loopref=\z@
     \loop\ifnum\@@loopref<\referencenum
               \advance\@@loopref by\@ne
               \setbox\z@=\hbox{\referencenum=\@@loopref
                    \referencenumformat\enskip}%
               \ifdim\wd\z@>\refnumindent\refnumindent=\wd\z@\fi
     \repeat
     \putreferenceformat
     \@@loopref=\z@
     \loop\ifnum\@@loopref<\@@prerefcount
               \advance\@@loopref by\@ne
               \expandafter\expandafter\expandafter\@printreference
                    \csname ref@\the\@@loopref\endcsname
     \repeat
     \let\reference=\@reference}

\def\@printreference#1#2{\ifx#2\undefinedlabelformat\err@undefinedref{#1}\fi
     \noindent\ifdraft\rlap{\hskip\hsize\hskip.1in \tt#1}\fi
     \llap{\referencenum=\@@loopref\referencenumformat\enskip}#2\par}

\def\reference#1#2{{\par\refnumindent=\z@\putreferenceformat\noindent#2\par}}

\def\putref#1{\trim#1\to\@@refarg
     \expandafter\ifnoarg\@@refarg\then
          \toks@={\relax}%
     \else\@@lastref=-\@m\def\@@refsep{}\def\@more{\@nextref}%
          \toks@={\@nextref#1,,}%
     \fi\the\toks@}
\def\@nextref#1,{\trim#1\to\@@refarg
     \expandafter\ifnoarg\@@refarg\then
          \let\@more=\relax
     \else\ifundefinedlabel\@@refarg\then
               \expandafter\@refpredef\expandafter{\@@refarg}%
          \fi
          \def\^^\##1{\global\@@thisref=##1\relax}%
          \global\@@thisref=\m@ne
          \setbox\z@=\hbox{\putlab\@@refarg}%
     \fi
     \advance\@@lastref by\@ne
     \ifnum\@@lastref=\@@thisref\advance\@@refseq by\@ne\else\@@refseq=\@ne\fi
     \ifnum\@@lastref<\z@
     \else\ifnum\@@refseq<\thr@@
               \@@refsep\def\@@refsep{,}%
               \ifnum\@@lastref>\z@
                    \advance\@@lastref by\m@ne
                    {\referencenum=\@@lastref\putrefformat}%
               \else\undefinedlabelformat
               \fi
          \else\def\@@refsep{--}%
          \fi
     \fi
     \@@lastref=\@@thisref
     \@more}

\message{streaming,}

\newif\ifstreaming

\def\streamto{\defaultoption[\jobname]\@streamto}
\def\@streamto[#1]{\global\streamingtrue
     \immediate\write\sixt@@n{--> Streaming to #1.str}%
     \newwrite\streamout\immediate\openout\streamout=#1.str }

\def\streamfrom{\defaultoption[\jobname]\@streamfrom}
\def\@streamfrom[#1]{\newread\streamin\openin\streamin=#1.str
     \ifeof\streamin
          \expandafter\err@nostream\expandafter{#1.str}%
     \else\immediate\write\sixt@@n{--> Streaming from #1.str}%
          \let\@@labeldef=\gdef
          \ifstreaming
               \edef\@elc{\endlinechar=\the\endlinechar}%
               \endlinechar=\m@ne
               \loop\read\streamin to\@@scratcha
                    \ifeof\streamin
                         \streamingfalse
                    \else\toks@=\expandafter{\@@scratcha}%
                         \immediate\write\streamout{\the\toks@}%
                    \fi
                    \ifstreaming
               \repeat
               \@elc
               \input #1.str
               \streamingtrue
          \else\input #1.str
          \fi
          \let\@@labeldef=\xdef
     \fi}

\def\streamcheck{\ifstreaming
     \immediate\write\streamout{\pagenum=\the\pagenum}%
     \immediate\write\streamout{\footnotenum=\the\footnotenum}%
     \immediate\write\streamout{\referencenum=\the\referencenum}%
     \immediate\write\streamout{\chapternum=\the\chapternum}%
     \immediate\write\streamout{\sectionnum=\the\sectionnum}%
     \immediate\write\streamout{\subsectionnum=\the\subsectionnum}%
     \immediate\write\streamout{\equationnum=\the\equationnum}%
     \immediate\write\streamout{\subequationnum=\the\subequationnum}%
     \immediate\write\streamout{\figurenum=\the\figurenum}%
     \immediate\write\streamout{\subfigurenum=\the\subfigurenum}%
     \immediate\write\streamout{\tablenum=\the\tablenum}%
     \immediate\write\streamout{\subtablenum=\the\subtablenum}%
     \immediate\closeout\streamout
     \fi}


\def\err@badtypesize{%
     \errhelp={The limited availability of certain fonts requires^^J%
          that the base type size be 10pt, 12pt, or 14pt.^^J}%
     \errmessage{--> Illegal base type size}}

\def\err@badsizechange{\immediate\write\sixt@@n
     {--> Size change not allowed in math mode, ignored}}

\def\err@sizetoolarge#1{\immediate\write\sixt@@n
     {--> \noexpand#1 too big, substituting HUGE}}

\def\err@sizenotavailable#1{\immediate\write\sixt@@n
     {--> Size not available, \noexpand#1 ignored}}

\def\err@fontnotavailable#1{\immediate\write\sixt@@n
     {--> Font not available, \noexpand#1 ignored}}

\def\err@sltoit{\immediate\write\sixt@@n
     {--> Style \noexpand\sl not available, substituting \noexpand\it}%
     \it}

\def\err@bfstobf{\immediate\write\sixt@@n
     {--> Style \noexpand\bfs not available, substituting \noexpand\bf}%
     \bf}

\def\err@badgroup#1#2{%
     \errhelp={The block you have just tried to close was not the one^^J%
          most recently opened.^^J}%
     \errmessage{--> \noexpand\end{#1} doesn't match \noexpand\begin{#2}}}

\def\err@badcountervalue#1{\immediate\write\sixt@@n
     {--> Counter (#1) out of bounds}}

\def\err@extrafootnotemark{\immediate\write\sixt@@n
     {--> \noexpand\footnotemark command
          has no corresponding \noexpand\footnotetext}}

\def\err@extrafootnotetext{%
     \errhelp{You have given a \noexpand\footnotetext command without first
          specifying^^Ja \noexpand\footnotemark.^^J}%
     \errmessage{--> \noexpand\footnotetext command has no corresponding
          \noexpand\footnotemark}}

\def\err@labelredef#1{\immediate\write\sixt@@n
     {--> Label "#1" redefined}}

\def\err@badlabelmatch#1{\immediate\write\sixt@@n
     {--> Definition of label "#1" doesn't match value in \jobname.lab}}

\def\err@needlabel#1{\immediate\write\sixt@@n
     {--> Label "#1" cited before its definition}}

\def\err@undefinedlabel#1{\immediate\write\sixt@@n
     {--> Label "#1" cited but never defined}}

\def\err@undefinedeqn#1{\immediate\write\sixt@@n
     {--> Equation "#1" not defined}}

\def\err@undefinedref#1{\immediate\write\sixt@@n
     {--> Reference "#1" not defined}}

\def\err@nostream#1{%
     \errhelp={You have tried to input a stream file that doesn't exist.^^J}%
     \errmessage{--> Stream file #1 not found}}

\message{jyTeX initialization}

\everyjob{\immediate\write16{--> jyTeX version \fmtversion}%
     \edef\@@jobname{\jobname}%
     \edef\jobname{\@@jobname}%
     \settime
     \openin0=\jobname.lab
     \ifeof0
     \else\closein0
          \immediate\write16{--> Getting labels from file \jobname.lab}%
          \input\jobname.lab
     \fi}


\def\fixedskipslist{%
     \^^\{\topskip}%
     \^^\{\splittopskip}%
     \^^\{\maxdepth}%
     \^^\{\skip\topins}%
     \^^\{\skip\footins}%
     \^^\{\headskip}%
     \^^\{\footskip}}

\def\scalingskipslist{%
     \^^\{\p@renwd}%
     \^^\{\delimitershortfall}%
     \^^\{\nulldelimiterspace}%
     \^^\{\scriptspace}%
     \^^\{\jot}%
     \^^\{\normalbaselineskip}%
     \^^\{\normallineskip}%
     \^^\{\normallineskiplimit}%
     \^^\{\baselineskip}%
     \^^\{\lineskip}%
     \^^\{\lineskiplimit}%
     \^^\{\bigskipamount}%
     \^^\{\medskipamount}%
     \^^\{\smallskipamount}%
     \^^\{\parskip}%
     \^^\{\parindent}%
     \^^\{\abovedisplayskip}%
     \^^\{\belowdisplayskip}%
     \^^\{\abovedisplayshortskip}%
     \^^\{\belowdisplayshortskip}%
     \^^\{\abovechapterskip}%
     \^^\{\belowchapterskip}%
     \^^\{\abovesectionskip}%
     \^^\{\belowsectionskip}%
     \^^\{\abovesubsectionskip}%
     \^^\{\belowsubsectionskip}}


\def\twoupsetup{
     \topmargin=.75in
     \leftmargin=.5in
     \vsize=6.9in
     \hsize=4.75in
     \fullhsize=10in
     \let\draft=\relax}

\outputstyle{normal}                             

\def\marginnoteformat{\subscriptsize             
     \hsize=1in \baselinestretch=1000 \everypar={}%
     \tolerance=5000 \hbadness=5000 \parskip=0pt \parindent=0pt
     \leftskip=0pt \rightskip=0pt \raggedright}

\head={\ifdraft\normalfonts\it\hfil DRAFT\hfil   
     \llap{\number\day\ \monthword\month\ \militarytime}\else\hfil\fi}
\foot={\hfil\normalfonts\numstyle\pagenum\hfil}  

\normalbaselineskip=12pt                         
\normallineskip=0pt                              
\normallineskiplimit=0pt                         
\normalbaselines                                 

\topskip=.85\baselineskip \splittopskip=\topskip \headskip=2\baselineskip
\footskip=\headskip

\pagenumstyle{arabic}                            

\parskip=0pt                                     
\parindent=20pt                                  

\baselinestretch=1000                            


\chapterstyle{left}                              
\chapternumstyle{blank}                          
\def\chapterbreak{\newpage}                      
\abovechapterskip=0pt                            
\belowchapterskip=1.5\baselineskip               
     plus.38\baselineskip minus.38\baselineskip
\def\chapternumformat{\numstyle\chapternum.}     

\sectionstyle{left}                              
\sectionnumstyle{blank}                          
\def\sectionbreak{\vskip0pt plus4\baselineskip\penalty-100
     \vskip0pt plus-4\baselineskip}              
\abovesectionskip=1.5\baselineskip               
     plus.38\baselineskip minus.38\baselineskip
\belowsectionskip=\the\baselineskip              
     plus.25\baselineskip minus.25\baselineskip
\def\sectionnumformat{
     \ifblank\chapternumstyle\then\else\numstyle\chapternum.\fi
     \numstyle\sectionnum.}

\subsectionstyle{left}                           
\subsectionnumstyle{blank}                       
\def\subsectionbreak{\vskip0pt plus4\baselineskip\penalty-100
     \vskip0pt plus-4\baselineskip}              
\abovesubsectionskip=\the\baselineskip           
     plus.25\baselineskip minus.25\baselineskip
\belowsubsectionskip=.75\baselineskip            
     plus.19\baselineskip minus.19\baselineskip
\def\subsectionnumformat{
     \ifblank\chapternumstyle\then\else\numstyle\chapternum.\fi
     \ifblank\sectionnumstyle\then\else\numstyle\sectionnum.\fi
     \numstyle\subsectionnum.}


\footnotenumstyle{symbol}                       
\footnoteskip=0pt                                
\def\footnotenumformat{\numstyle\footnotenum}    
\def\footnoteformat{\footnotesize                
     \everypar={}\parskip=0pt \parfillskip=0pt plus1fil
     \leftskip=1em \rightskip=0pt
     \spaceskip=0pt \xspaceskip=0pt
     \def\\{\ifhmode\ifnum\lastpenalty=-10000
          \else\hfil\penalty-10000 \fi\fi\ignorespaces}}


\def\undefinedlabelformat{$\bullet$}             


\equationnumstyle{arabic}                        
\subequationnumstyle{blank}                      
\figurenumstyle{arabic}                          
\subfigurenumstyle{blank}                        
\tablenumstyle{arabic}                           
\subtablenumstyle{blank}                         

\eqnseriesstyle{alphabetic}                      
\figseriesstyle{alphabetic}                      
\tblseriesstyle{alphabetic}                      

\def\puteqnformat{\hbox{
     \ifblank\chapternumstyle\then\else\numstyle\chapternum.\fi
     \ifblank\sectionnumstyle\then\else\numstyle\sectionnum.\fi
     \ifblank\subsectionnumstyle\then\else\numstyle\subsectionnum.\fi
     \numstyle\equationnum
     \numstyle\subequationnum}}
\def\putfigformat{\hbox{
     \ifblank\chapternumstyle\then\else\numstyle\chapternum.\fi
     \ifblank\sectionnumstyle\then\else\numstyle\sectionnum.\fi
     \ifblank\subsectionnumstyle\then\else\numstyle\subsectionnum.\fi
     \numstyle\figurenum
     \numstyle\subfigurenum}}
\def\puttblformat{\hbox{
     \ifblank\chapternumstyle\then\else\numstyle\chapternum.\fi
     \ifblank\sectionnumstyle\then\else\numstyle\sectionnum.\fi
     \ifblank\subsectionnumstyle\then\else\numstyle\subsectionnum.\fi
     \numstyle\tablenum
     \numstyle\subtablenum}}


\referencestyle{sequential}                      
\referencenumstyle{arabic}                       
\def\putrefformat{\numstyle\referencenum}        
\def\referencenumformat{\numstyle\referencenum.} 
\def\putreferenceformat{
     \everypar={\hangindent=1em \hangafter=1 }%
     \def\\{\hfil\break\null\hskip-1em \ignorespaces}%
     \leftskip=\refnumindent\parindent=0pt \interlinepenalty=1000 }


\normalsize


\def\fmtversion{2.6M (June 1992)}

\catcode`\@=12
\typesize=10pt \magnification=1200 \baselineskip17truept
\footnotenumstyle{arabic} \hsize=6truein\vsize=8.5truein
\sectionnumstyle{blank}
\chapternumstyle{blank}
\chapternum=1
\sectionnum=1
\pagenum=0

\def\begintitle{\pagenumstyle{blank}\parindent=0pt
\begin{narrow}[0.4in]}
\def\endtitle{\end{narrow}\newpage\pagenumstyle{arabic}}


\def\beginexercise{\vskip 20truept\parindent=0pt\begin{narrow}[10
truept]}
\def\endexercise{\vskip 10truept\end{narrow}}


\def\eql#1{\eqno\eqnlabel{#1}}
\def\ref{\reference}
\def\peq{\puteqn}
\def\pref{\putref}

\def\mgn{\marginnote}
\def\bex{\begin{exercise}}
\def\eex{\end{exercise}}


\font\open=msbm10 


\def\StretchRtArr#1{{\count255=0\loop\relbar\joinrel\advance\count255 by1
\ifnum\count255<#1\repeat\rightarrow}}
\def\StretchLtArr#1{\,{\leftarrow\!\!\count255=0\loop\relbar
\joinrel\advance\count255 by1\ifnum\count255<#1\repeat}}

\def\StretchLRtArr#1{\,{\leftarrow\!\!\count255=0\loop\relbar\joinrel\advance
\count255 by1\ifnum\count255<#1\repeat\rightarrow\,\,}}

\def\mbox#1{{\leavevmode\hbox{#1}}}

\def\hspace#1{{\phantom{\mbox#1}}}
\def\oR{\mbox{\open\char82}}

\def\oZ{\mbox{\open\char90}}

\def\al{\alpha}
\def\bde{{\bmit\delta}} 
\def\bDe{{\bmit\Delta}}
\def\be{\beta}

\def\de{\delta}
\def\Ga{\Gamma}

\def\ep{\epsilon}

\def\la{\lambda}

\def\th{\theta}

\def\De{\Delta}

\def\caN{{\cal N}}

\def\caD{{\cal D}}

\def\sc{{\rm sc }}


\def\frac#1/#2{\leavevmode\kern.1em
\raise.5ex\hbox{\the\scriptfont0 #1}\kern-.1em/\kern-.15em
\lower.25ex\hbox{\the\scriptfont0 #2}}
\def\sfrac#1/#2{\leavevmode\kern.1em
\raise.5ex\hbox{\the\scriptscriptfont0 #1}\kern-.1em/\kern-.15em
\lower.25ex\hbox{\the\scriptscriptfont0 #2}}

\def\gtorder{\mathrel{\raise.3ex\hbox{$>$}\mkern-14mu
             \lower0.6ex\hbox{$\sim$}}}
\def\ltorder{\mathrel{\raise.3ex\hbox{$<$}\mkern-14mu
             \lower0.6ex\hbox{$\sim$}}}

\def\semidirprod{\rlap{\ss C}\raise1pt\hbox{$\mkern.75mu\times$}}
\def\for{\lower6pt\hbox{$\Big|$}}
\def\fish{\kern-.25em{\phantom{abcde}\over \phantom{abcde}}\kern-.25em}


\def\boxit#1{\vbox{\hrule\hbox{\vrule\kern3pt
        \vbox{\kern3pt#1\kern3pt}\kern3pt\vrule}\hrule}}
\def\dalemb#1#2{{\vbox{\hrule height .#2pt
        \hbox{\vrule width.#2pt height#1pt \kern#1pt \vrule
                width.#2pt} \hrule height.#2pt}}}

\def\ol{\overline}
\def\frac#1#2{{{#1}\over{#2}}}

\def\noin{\noindent}


\def\eg{{\it e.g.}}
\def\ie{{\it i.e. }}
\def\cf{{\it cf }}
\def\pa{\partial}

\def\ket#1{\mid#1\rangle}
\def\br#1#2{\langle{#1}\mid{#2}\rangle}   

\def\wt{\widetilde}

\def\curl{{\rm curl\,}}
\def\div{{\rm div\,}}
\def\grad{{\rm grad\,}}

\def\Threej#1#2#3#4#5#6{\biggl({#1\atop#4}{#2\atop#5}{#3\atop#6}\biggr)}
\def\threej#1#2#3#4#5#6{\bigl({#1\atop#4}{#2\atop#5}{#3\atop#6}\bigr)}
\def\3j#1#2#3#4#5#6{\left\lgroup\matrix{#1&#2&#3\cr#4&#5&#6\cr}
\right\rgroup}
\def\sixj#1#2#3#4#5#6{\left\{\matrix{\!#1&\!\!
\!#2&\!\!\!#3\!\cr\!#4&\!\!\!#5&\!\!\!#6\!\cr}\right\}}

\def\6j#1#2#3#4#5#6{\left\{{#1\atop#4}{#2\atop#5}{#3\atop#6}\right\}}

\def\m?{\mgn{?}}

\def\pa{\partial}

\def\beq{\begin{eqnarray}}
\def\eeq{\end{eqnarray}}


\def\aop#1#2#3{{\it Ann. Phys.} {\bf {#1}} ({#2}) #3}

\def\cmp#1#2#3{{\it Comm. Math. Phys.} {\bf {#1}} ({#2}) #3}
\def\cqg#1#2#3{{\it Class. Quant. Grav.} {\bf {#1}} ({#2}) #3}

\def\ijmp#1#2#3{{\it Int. J. Mod. Phys.} {\bf {#1}} ({#2}) #3}

\def\jmp#1#2#3{{\it J. Math. Phys.} {\bf {#1}} ({#2}) #3}
\def\jpa#1#2#3{{\it J. Phys.} {\bf A{#1}} ({#2}) #3}
\def\lnm#1#2#3{{\it Lect. Notes Math.} {\bf {#1}} ({#2}) #3}

\def\np#1#2#3{{\it Nucl. Phys.} {\bf B{#1}} ({#2}) #3}
\def\pl#1#2#3{{\it Phys. Lett.} {\bf {#1}} ({#2}) #3}

\def\prp#1#2#3{{\it Phys. Rep.} {\bf {#1}} ({#2}) #3}
\def\pr#1#2#3{{\it Phys. Rev.} {\bf {#1}} ({#2}) #3}
\def\prA#1#2#3{{\it Phys. Rev.} {\bf A{#1}} ({#2}) #3}

\def\prD#1#2#3{{\it Phys. Rev.} {\bf D{#1}} ({#2}) #3}
\def\prl#1#2#3{{\it Phys. Rev. Lett.} {\bf #1} ({#2}) #3}

\def\rmp#1#2#3{{\it Rev. Mod. Phys.} {\bf {#1}} ({#2}) #3}

\def\zfp#1#2#3{{\it Z. f. Phys.} {\bf {#1}} ({#2}) #3}

\def\cras#1#2#3{{\it Comptes Rend. Acad. Sci. (Paris)} {\bf{#1}} (#2) #3}
\def\prs#1#2#3{{\it Proc. Roy. Soc.} {\bf A{#1}} ({#2}) #3}
\def\pcps#1#2#3{{\it Proc. Camb. Phil. Soc.} {\bf{#1}} ({#2}) #3}
\def\mpcps#1#2#3{{\it Math. Proc. Camb. Phil. Soc.} {\bf{#1}} ({#2}) #3}

\def\amsh#1#2#3{{\it Abh. Math. Sem. Ham.} {\bf {#1}} ({#2}) #3}
\def\am#1#2#3{{\it Acta Mathematica} {\bf {#1}} ({#2}) #3}
\def\aim#1#2#3{{\it Adv. in Math.} {\bf {#1}} ({#2}) #3}
\def\ajm#1#2#3{{\it Am. J. Math.} {\bf {#1}} ({#2}) #3}

\def\aom#1#2#3{{\it Ann. of Math.} {\bf {#1}} ({#2}) #3}
\def\cjm#1#2#3{{\it Can. J. Math.} {\bf {#1}} ({#2}) #3}
\def\bams#1#2#3{{\it Bull.Am.Math.Soc.} {\bf {#1}} ({#2}) #3}

\def\cmh#1#2#3{{\it Comm. Math. Helv.} {\bf {#1}} ({#2}) #3}

\def\dmj#1#2#3{{\it Duke Math. J.} {\bf {#1}} ({#2}) #3}
\def\invm#1#2#3{{\it Invent. Math.} {\bf {#1}} ({#2}) #3}

\def\jdg#1#2#3{{\it J. Diff. Geom.} {\bf {#1}} ({#2}) #3}

\def\joa#1#2#3{{\it J. of Algebra} {\bf {#1}} ({#2}) #3}
\def\jram#1#2#3{{\it J. f. reine u. Angew. Math.} {\bf {#1}} ({#2}) #3}
\def\jims#1#2#3{{\it J. Indian. Math. Soc.} {\bf {#1}} ({#2}) #3}
\def\jlms#1#2#3{{\it J. Lond. Math. Soc.} {\bf {#1}} ({#2}) #3}
\def\jmpa#1#2#3{{\it J. Math. Pures. Appl.} {\bf {#1}} ({#2}) #3}
\def\ma#1#2#3{{\it Math. Ann.} {\bf {#1}} ({#2}) #3}

\def\mz#1#2#3{{\it Math. Zeit.} {\bf {#1}} ({#2}) #3}
\def\ojm#1#2#3{{\it Osaka J.Math.} {\bf {#1}} ({#2}) #3}

\def\pems#1#2#3{{\it Proc. Edin. Math. Soc.} {\bf {#1}} ({#2}) #3}

\def\plb#1#2#3{{\it Phys. Letts.} {\bf {B#1}} ({#2}) #3}
\def\pla#1#2#3{{\it Phys. Letts.} {\bf {A#1}} ({#2}) #3}
\def\plms#1#2#3{{\it Proc. Lond. Math. Soc.} {\bf {#1}} ({#2}) #3}
\def\pgma#1#2#3{{\it Proc. Glasgow Math. Ass.} {\bf {#1}} ({#2}) #3}
\def\qjm#1#2#3{{\it Quart. J. Math.} {\bf {#1}} ({#2}) #3}
\def\qjpam#1#2#3{{\it Quart. J. Pure and Appl. Math.} {\bf {#1}} ({#2}) #3}

\def\rmjm#1#2#3{{\it Rocky Mountain J. Math.} {\bf {#1}} ({#2}) #3}

\def\tams#1#2#3{{\it Trans.Am.Math.Soc.} {\bf {#1}} ({#2}) #3}

\input epsf

\begin{title}
\vglue 1truein
\vskip15truept
\centertext {\Bigfonts \bf Algebraic derivation of some theorems concerning} \vskip10truept
\centertext{\Bigfonts \bf contact structures on the 3--sphere }
 \vskip 20truept
\centertext{J.S.Dowker\footnote{dowkeruk@yahoo.co.uk}} \vskip 7truept \centertext{\it
Theory Group,} \centertext{\it Department of Physics and Astronomy,} \centertext{\it The
University of Manchester,} \centertext{\it Manchester, England} \vskip40truept
\begin{narrow}
Two  theorems involving curl eigenfields on the 3--sphere are obtained using angular momentum theory. Spinor hyperspherical harmonics are shown to  form an explicit, convenient basis. In particular,  a spin--one vector calculus is reviewed. An easy proof of the vanishing of  `odd' eigenfields is given and related to the sign change of fermionic spinors under 2$\pi$ rotations. The theorem that curl eigenfields with constant norm have to be proportional to a fundamental eigenfield (Hopf field) is also rapidly obtained.

Attention is drawn to the relevance of early work of Schr\"odinger on Maxwell theory in an expanding universe.

\end{narrow}

\vskip 5truept
\vskip 60truept
\vfil
\end{title}
\pagenum=0
\newpage
\section{\bf1. Introduction}

In [\pref{DandP2}], in order to compute the Schr\"odinger energy eigenvalues of an ideal asymmetric top (with spin) or, equivalently, of a higher--spin field in a Mixmaster universe, spinor hyperspherical harmonics on the three--sphere, S$^3$, were introduced. These can be considered to be generalisations both of the (scalar) Wigner curly $\caD$ functions on S$^3$ and of the standard spinor spherical harmonics on S$^2$ (\eg\ [\pref{LandL}] [\pref{Lyubarskii}]). They were employed, [\pref{Dowk}], in a discussion of arbitrary spin fields on the Einstein universe T$\times$ S$^3$. For massless fields, for example Maxwell (spin--one), the time development operator takes the form of a higher--spin curl whose eigenmodes were given explicitly by spinor hyperspherical harmonics. Further, the spin--one curl eigenfields, in particular, were used, [\pref{Dowk3,Dowk5}], as a convenient basis for the diffeomorphisms and fluid flow on S$^3$. The computation of invariant quantities was thereby efficiently reduced to one in Wigner--Racah calculus.

Relatedly, curl eigenfields appear, somewhat fundamentally, in the theory of contact structures and hydrodynamics. In dimension three, the only manifolds (of positive curvature) having a normal contact structure are diffeomorphic to S$^3/\Gamma$ with a freely acting $\Gamma$, [\pref{Geiges}]. 

The paper [\pref{PSS}] is concerned, partly, with the standard round sphere  ($\Gamma={\bf 1}$) and proves several theorems which depend on a particular product construction of the higher--eigenvalue curl eigenfields. The aim of the present short note is to re--examine some of the results and methods ln [\pref{PSS}] using, instead, the spinor hyperspherical basis and basic angular momentum theory.

The next three sections review the higher spin modes and some of their properties. The novel application to contact structures is contained in the shortish section 5.

\section{\bf 2. Modes on the three--sphere}

Historically, the Maxwell eigenfunction problem on the Einstein universe, T$\times$S$^3$, was completely solved by Schr\"odinger\footnote{ With help from Pauli}, [\pref{Schr}], in 1940 in the traditional way of selecting a coordinate system and solving the associated partial differential equation by separation of variables and application of regularity. \footnote{ In an earlier paper, [\pref{Schr2}], he treats scalar waves by group theory using the SU(2)$\times$SU)2)$\sim$ SO(4) symmetry.}  (Actually he treats the case when the radius of the sphere is time dependent, but shows that this is not an obstacle to a solution by redefining time to $\tau$. This does not concern me here.) A cylindrical coordinate system\footnote{ \cf\ Wenger, [\pref{Wenger}]. This is the Hopf coordinate system also used in [\pref{PSS}]} is chosen and  the curl eigenvalue computed in the guise of the frequency, $\nu$, this being the eigenvalue of $-i\pa/\pa\tau$. After a lengthy analysis (as Schr\"odinger admits), $|\nu|$ is shown to be an integer greater than or equal to 2 and a rule is given for finding the multiplicity which turns out to be $\nu^2-1$. He notes the similarity of this with $\nu^2$ for spin 0 and $\nu^2-(1/2)^2$ for spin-1/2.

Schr\"odinger gives an extended description of the fundamental mode \ie $\nu=2$, and notes that the $\nu=-2$ mode is just its mirror image. This discussion is very pertinant\footnote{ It is interesting to note that Schr\"odinger, in App.V, produces solutions of Maxwell's eqations (but not curl eigenfunctions) which are obtained by multiplying the fundamental one by an arbitrary analytic function, a process akin to that employed in [\pref{PSS}].} since, when normalised, this mode, choosing one of the three possible,  is just the Reeb vector field of the standard contact structure on S$^3$. In classical Lie group theory terms, this field is one of the left--invariant, Killing vector fields, $A^\al_a$, labelled by $a=(1,0,-1)$, say and where $\al$ is a general coordinate index on S$^3\sim$SU(2). It is easily checked directly from the Maurer-Cartan equation that this field is divergence--free and has a curl eigenvalue of 2. Supplied with a derivative it gives the generators of the right regular representation of SU(2), $Y_a\equiv A_a^\al\pa_\al$, which can be taken as the Reeb vector field(s).
A parity inversion takes us to the left generators with  the corresponding vector fields having curl eigenvalue, $ -2$.

All this is fairly basic and, to further the analysis, I now introduce the notion of right/left $j$--spinors, $\phi^j_m$, which are such that under the left--right symmetry group action $q\to q'=\xi q \eta$ where $q$ and $q'$ belong to SU(2), \ie points on S$^3$, the $\phi$ transform according to,\footnote{ The $\caD$ are the usual Wigner functions. The second projection number has been written upstairs purely for cosmetic, index tracking purposes. It is not raised by any metric.}
  $$\eqalign{
  \phi^j_m(q)&\to\phi'_{m'}(q')={\caD^J_{m'}}^{m}(\eta^{-1})\,\phi_{m}(q)\,,\quad
  {\rm right}\cr
  \phi^j_m(q)&\to\phi'_{m'}(q')={\caD^J_{m'}}^{m}(\xi)\,\phi_{m}(q)\,,\quad
  {\rm left}\,.\cr
  }
  \eql{lr}
  $$

For reasons described below, it is algebraically more convenient to choose right fields,\footnote{ Of course there are many, ultimately equivalent ways of defining spinors because, in general, spin--space is completely independent of coordinate space.} a complete set of modes for which is provided by the right spinor hyperspherical harmonics. These have the coordinate representation, [\pref{DandP2}], [\pref{Pett}],
$$\eqalign{
  Y^{j\,LM}_{mNJ}(q)&=\br{m,q}{JLj(N)M}\cr
  &=\bigg({(2L+1)(2J+1)\over2\pi^2R^3}\bigg)^{1/2}\caD^{LN'}_N(q)
  \Threej jLMm{N'}J\,,
  }
  \eql{hsh}
	$$
which is a straightforward generalisation  of the usual spin--orbit coupling with ${\bf j}$ being the spin, ${\bf L}$ the right "orbital" angular momentum (${\bf L}=-i{\bf Y}$) and ${\bf J},={\bf L}+{\bf j}$, the total angular momentum. The states $\ket{JLj(N)M}$ have degeneracy $(2L+1)(2J+1)$ and  diagonalise the operator, 
$$
   H(a,b,c)=a{\bf
   L}^2+b{\bf J}^2 +c{\bf j}^2\,,
   \eql{lap4}
   $$
which can be realised as the Hamiltonian of an ideal spherical top, with spin, [\pref{DandP2}]. or as a higher--spin curl operator or as  a Laplacian (see section (3).

A technical point is that, when the radius, $R$, of S$^3$ equals 2, the states\break $\ket{JLj(N)M}$ can be used in all the usual angular momentum calculations with no changes in conventions, or normalisation.\footnote{ Left spinors can be defined similarly but then a complex conjugation step is needed.} The left label $N$ is simply a spectator as indicated by the bracket. I therefore set the radius to 2 which accounts for any differences in factors of 2 = $R$ and 4 = $R^2$ with other works which usually assume a unit S$^3$. For example, here, the Reeb vector curl eigenvalues are $\pm 1$.

For spin 1, the right--spinor is just the projection of the vector field onto the left invariant Killing vector field,
  $$
	   V_a=A^\al_a V_\al\,.
		\eql{proj}
	$$
and the inverse,
  $$
		V_\al=A^a_\al V_a\,,\quad {\rm where}\quad A^a_\al\,A_a^\be=\de_\be^\al\,,\quad{\rm or}\quad
		A^a_\al A^b_\be\,g_{ab}=g_{\al\be}\,.
		\eql{projinv}
	$$	
Here, $g_{\al\be}$ is the standard metric on S$^3$ and $g_{ab}$ the Cartan metric on adjoint, here spin-one, space.	\footnote{ This last can be taken to be the Euclidean unit matrix or the angular momentum spherical one.}

The projection onto the right--invariant	frame is,
$$
	   \wt V_a=\wt A^\al_a V_\al\,.
		\eql{projr}
	$$

Left and right are related by an adjoint (spin--one) rotation,
  $$
	  \wt A^a_\al(q)=\caD^{1a}_b(q)\,A^b_\al(q)\,.
	$$

\section{\bf 3. The operators Laplacian and curl}

To construct covariant or invariant equations of motion, a covariant derivative is required. For right spinors this was found in [\pref{Dowk2}] to be,
  $$
	{\bf\nabla}={\bf Y}+{1\over2} i\, {\bf j}\,,
	$$
so that a natural spinor Laplacian is,
   $$
	      \bDe_2\equiv {\bf\nabla}^2=-({\bf L}+{1\over2}{\bf j})^2=-{1\over4}\big( {2\bf L}^2+2{\bf J}^2-{\bf j}^2\big)\,,
				\eql{splap}
	$$
which is of the form of (\peq{lap4}) with eigenvalues $-J(J+1)/2-L(L+1)/2+j(j+1)/4$ so long as the triangle relation $\De(J,L,j)$ for angular momentum addition holds.

In order to motivate the expression for the curl operator, consider massless field equations in a curved space--time. There are numerous higher--spin generalizations of Maxwell theory in flat space--time, some are given in [\pref{DandD2}]. One which extends easily into curved space--time by minimal coupling is due to Dirac and then reads,
   $$
	 (iJ^{\mu\nu}\nabla_\nu+j\nabla^\mu)\phi=0\,,
	\eql{md}
	$$
where $\phi$ is a spin--j field belonging to the $(j,0)$ representation of the homogeneous (local) Lorentz group and $J^{\mu\nu}$ are the corresponding, self--dual generators. 

This equation is inconsistent in a generally curved space--time, except for spin--one. However the inconsistencies disappear if space--time is conformally flat and so can be safely applied to the Einstein universe\footnote{ and de Sitter space}, where it becomes, [\pref{Dowk}], using the left local basis,
     $$
		      \big (i\pa_0 +{1\over j} {\bf j}\cdot{\bf \nabla}\big)\phi=0
					\eql{curl}
		$$
and
    $$					
					\Threej a {j-1} {m'}1mj {\bf \nabla}_{a}\phi=0\,.
					\eql{div}
		$$
The equation has been split into a "curl" part, (\peq{curl}), and a "div" part, (\peq{div}). The latter,  on expanding $\phi$ in spinor hyperspherical harmonics, reduces to the triangle relation, $\De(J,L,j-1)$ which implies the massless polarisation conditions,
  $$
	J=L\pm j,\quad {\rm for}\,\, J\le L\,\quad 
	{\rm and}\quad J=L+j,\quad {\rm for} \,\,J>L\,.
	\eql{poln}
	$$

Although I am interested most in spin--one, the analysis proceeds easily for any $j$. The energies (frequencies) can be determined from (\peq{curl}) as the eigenvalues of $i\pa_0$ \ie\ of the generalised curl operator which is,
  $$
\curl\equiv	-i{1\over j} {\bf j}\cdot{\bf \nabla}={1\over j}\big({\bf j}\cdot{\bf L}+{1\over2}{\bf j}^2\big)\,,
	$$
with eigenvalues,
$$
       E_{L,J}^j={1\over 2j} (J-L)(J+L+1)  
			\eql{ceig}
$$
subject to the polarisation conditions, (\peq{poln}), which lead to the antisymmetrical energy spectrum\footnote{ I have reinstated, briefly, the radius $R$. Note that these values hold for spin zero although the equation of motion, (\peq{md}), makes no sense. However, the second order operator obtained by applying $\nabla_\mu$ to (\peq{md}) simplifies to $\De_2+R(j+1)/6$, the usual conformally covariant one for spin zero.}, [\pref{Pett}], [\pref{Dowk}],
  $$\eqalign{
	E^+_{\ol L}&=\pm{1\over R}(j+\ol L)\,\quad \ol L=1,2,3\ldots\,.\cr
	}
	\eql{spect}
	$$
The degeneracies are $(2L+1)(2J+1)$.

The modes with equal and opposite energies are obtained by interchanging $L$ and $J$ and are parity related under the reflection in the unit group element, \ie\ $q\to q^{-1}$, [\pref{Dowk}], I will, therefore, concentrate on just the positive modes.

\section{\bf4. Spin one calculus}

A convenient vector calculus can be set up for spin one, [\pref{Dowk3,Arnold}], that parallels the Gibbs--Heaviside calculus in flat space.\footnote{ A $p$--form description is equally effective.} I expand a little on the discussion in [\pref{Dowk3}].

The three basic operations of $\curl$, $\div$ and $\grad$ are defined on right 3--spinors, ${\bf h}$, and scalars, $\phi$, by,
  $$
  \curl {\bf h}={\bf Y}\times{\bf h}+{\bf h}\,,\quad \div {\bf h}={\bf
Y\cdot h}\,,\quad \grad\phi={\bf Y}\phi\,,
  $$
and follow, for example, from their covariant Riemannian expressions restricted to the
three--sphere and projected onto the left--invariant Killing fields. The cross product is the usual one, the standard $Y$--commutator being  ${\bf Y\times Y}=-{\bf Y}$.

The usual identities hold, as is easily shown. Thus,
  $$
  \curl\grad\phi={\bf Y\times Y}\phi+{\bf Y}\phi=-{\bf Y}\phi+{\bf
Y}\phi=0\,,
  $$
and
  $$\eqalign{
  \div\curl{\bf h}&={\bf Y .Y}\times{\bf h}+{\bf Y.h}
  ={\bf Y \times Y}.{\bf h}+{\bf Y.h}=-{\bf Y.h}+{\bf Y.h}\cr
  &=0\,.
}
  $$

Also the repeated curl equation,
  $$\eqalign{
  \curl\curl{\bf h}&={\bf Y}\times({\bf Y\times h+h})+{\bf Y\times h+h}\cr
 &={\bf Y}({\bf Y.h})-{\bf Y\times h}-{\bf Y}^2{\bf h} +2{\bf Y\times
h+h}\cr
 &=\grad\div{\bf h}-\bigg[\big({\bf Y}+{i\over2}\,{\bf
 j}\big)^2-{1\over2}\bigg]{\bf h}\,,
} \eql{curl2}
  $$
generalises the standard relation on $\oR^3$, \cf [\pref{Arnold}].

A few manipulations have been made in order to introduce the spinor Laplacian, ({\peq{splap}).
In particular, I remark that the ${\bf j}$ are the spin--1 angular momentum matrices which
are the adjoint structure constants, $\ep$, up to a factor of $i$. Thus $i{\bf Y.\,j}\equiv {\bf
Y}\times$. The $-{\bf Y\times h}$ term arises from a commutator required to give the grad
div term. The expression in square brackets is the de Rham Laplacian, $\De_{\rm
deR}=d\de+\de d$ acting on 1--forms.

Any right vector field, ${\bf V}$, on S$^3$ can be expanded in spin-1 harmonics (modes) defined by\footnote{ I apologise for the dual use of the symbol {\bf Y}.},
  $$
  \big({\bf Y}^{L\,M}_{NJ}\big)_m\equiv Y^{1\,L\,M}_{mNJ},,
	\eql{vy}
  $$
as particular spinor hyperspherical harmonics, (\peq{hsh}). Thus,
  $$
  {\bf V}=\sum_{LNJM}V^{NJ}_{LM}\,{\bf Y}^{L\,M}_{NJ}\,.
  \eql{sp1exp}
  $$

Angular momentum addition implies ${\bf J=L+1}$. The modes with $J=L$ are, directly
from their definition, longitudinal, \ie gradients (=$\grad\phi$), by virtue of the basic right action equation, 
  $$
  -i{\bf Y}\caD^J(q)=\caD^J(q)\,{\bf J}\,,
  \eql{righta}
  $$
and the fact that the $3j$ symbols, $\threej 1LMmNL$, are
proportional to the angular momentum matrices in the $L$--representation (see below). They are zero modes of the curl operator, see (\peq{ceig}). The actual relation is
  $$
  {\bf Y}^{L\,M}_{N\,L}(q)={1\over\sqrt{L(L+1}}\,{\bf Y}\,Y^{0LM}_{\,0\,NL}(q)
  \eql{long}
  $$
which, as said, follows from the generators,
  $$
  ({\bf L}_m)_N^{\,M}=i\big(L(L+1)(2L+1)\big)^{1/2}\Threej L1{M}NmL\,,
  $$
the action,
  $
  {\bf Y}\caD^L=i\caD^L{\bf L}\,,
  $
and the spin--zero modes which are, from the definition (\peq{hsh}),
  $$
   Y^{0\,LM}_{\,0\,NJ}(q)=(-1)^{2L}\,\de_{LJ}\,
  {\sqrt{2L+1}\over4\pi}\,\caD^{\,LM}_N(q)
  \eql{scmodes}
  $$
where the (inessential) phase, $(-1)^{2L}$ arises from the $3j$ symbol,
  $$
  \Threej 0LM0{N'}J=(-1)^{2L}\,\sqrt{2L+1}\,\de_{LJ}\,\,\de_{N'}^M\,.
  $$
	
	Finally, the projection to give a vector field, (\peq{proj}), can be written,
	  $$
		{V}_\al= {\bf A}_\al\cdot{\bf V}\,.
		$$
	
\section{\bf 5. More about modes. Applications to contact geometry}

The statements in section 2 can now be recast and extended. For example, the curl fundamental mode, from (\peq{spect}), is when $L=0$. and the harmonic is constant  To spell things out,
  $$\eqalign{
	  Y^{1\,0\,M}_{\,a\,0\,J}(q)&=\bigg({3\over16\pi^2}\bigg)^{1/2}
		\Threej 10Ma0J\cr
		&={3\over4\pi}\,\de_{J1}\,\de_a^M\,.
		}
	$$
or as a vector,
  $$
	{\bf  Y}^{0\,M}_{\,0\,J}(q)={3\over4\pi}\,\de_{J1}\,\bde^M\,,
	$$
which is the standard unit spinor in the local frame, \eg\ [\pref{Lyubarskii}].

When projected according to (\peq{proj}), this yields, of course,  the three possible positive Reeb vector eigenfields, as the left--invariant Killing vectors, as expected,
  $$
	   {3\over4\pi}A_\al^M(q)\,,\quad M=1,0,-1\,.
	$$
\begin{ignore}
For completness, I discuss the corresponding negative mode. This is obtained by interchanging $L$ and $J$ and gives a mode with $L=1$ and $J=0$, \ie,
  $$\eqalign{
	      Y^{1\,1\,0}_{\,a\,N\,0}(q)&=\bigg({3\over16\pi^2}\bigg)^{1/2}\caD^{1N'}_N(q)
		\Threej 110a{N'}0\cr
		&={3\over4\pi}\,\de_a^{-N'}(-1)^{N'}\,\caD^{1,N'}_N(q)\,.
		}
	$$

Projecting onto the Killing field gives the eigenfield,
  $$
	     A_\al^a \de_a^{-N'}(-1)^{N'}\,\caD^{1N'}_N(q)\,.
	$$
\end{ignore}

In [\pref{PSS}], these  are referred to as Hopf fields and linear combinations of them, with function coefficients, are used to construct the {\it higher} curl eigenfields. On a 3--Sasakian manifold, the coefficients are shown, in [\pref{PSS2}], to be eigenfunctions of the (scalar) Laplacian using  a slightly involved differential geometric method. By contrast, on S$^3$ (and S$^3/\Ga$), this is built into the form, (\peq{hsh}), of the curl eigenmodes. Consequently, it needs no special proof. (See the Added Note.)

This product structure is then applied in [\pref{PSS}] to show that eigenfields with odd eigenvalues (in the normalisation of [\pref{PSS}]) are vanishing (somewhere). In the development here, the positive curl eigenvalues are $E=1+L$, $L=0,1/2,1\dots$, and the `odd' eigenforms correspond to $L$ half--integral. All such eigenforms (and their linear combinations) have a zero as a conseqence of the famous sign change of spin half--integral (fermionic) functions under the, $\oZ_2$, antipodal map on S$^3$ taking the identity to rotation by $2\pi$. This is equivalent to appealing to the Borsuk--Ulam theorem as in [\pref{PSS}].

Another result proved in [\pref{PSS}], section 4 (Proposition 1), is that if ${\bf V}$ has a constant (\ie position independent) norm, then it is proportional to the fundamental mode (Hopf field). The proof again involves their product form of the higher modes,  as decribed above, and employs a harmonic mapping from S$^3$ to S$^2$, and other ingredients, to show that these modes cannot have constant norm. Here, the explicit expansion, (\peq{sp1exp}), can be used to advantage and reduces the calculation to a purely algebraic one in that angular momentum theory allows one to calculate the norm of ${\bf V}$ quite easily as follows.

The norm of any ${\bf V}$ is determined by,
$$
      |{\bf V}|^2={\bf V}\cdot{\bf V}=V^aV_a=\sum_{LNJM}\,\sum_{L'N'J'M'}V^{NJ}_{LM}V^{N'J'}_{L'M'}\,{\bf Y}^{L\,M}_{NJ}\cdot{\bf Y}^{L'\,M'}_{N'J'}\,,
$$
where indices $a$ are raised and lowered by the Wigner metric (here in adjoint, spin--one, space) and so one needs the general mode scalar product,
  $$
	\caN^{LML'M'}_{NJN'J'}\equiv{\bf Y}^{L\,M}_{NJ}\cdot{\bf Y}^{L'\,M'}_{N'J'}\,,
	\eql{scprod}
	$$
which will follow from the definition, (\peq{hsh}), and the standard Clebsch--Gordan product for $\caD\caD$ \eg\ [\pref{BandS}] equn (2.32). This last contains a single $\caD$ and two 3j symbols which, because of the scalar product, can be re--coupled to give a sum over a 6j symbol multiplied by just one 3j--symbol. 

Initially, this calculation is best performed graphically. Once the structure of the answer is known, the algebra can then be done to get the correct phases, if needed, which here they are not.

The result is, after this short rearrangement,\footnote{ Sum over $k,k'$ implied.}
  $$
	\caN^{LML'M'}_{NJN'J'}(q)\!=\!\sum_{K={\rm Max}}^{\rm Min}(-1)^\Pi\,\sixj J{J'}K{L'}L1\Threej L{L'}kN{N'}K
	\caD^{Kk'}_k(q)\Threej KM{M'}{k'}J{J'}\,,
	\eql{norm1}
	$$
to an overall constant. The limits are Min = Min$(L+L',J+J')$ and Max  = Max$(|L-L'|,|J-J'|)$. $\Pi$ is a phase which we do not need unless a numerical answer is sought. $K$ is either integral or half--integral depending on the values  of $L$ and $L'$.

We are interested especially in the case when both modes in (\peq{scprod}) refer to the same, positive curl eigenvalue. That is, when $L'=L$ and, since $J=L+1$, also  when $J'=J$. The magnetic quantum numbers $N,M,N',M'$ can be left free so as to allow the construction of a general curl eigenmode if and when required. 

$K$ must now be integral and it is immediately apparent from (\peq{norm1}), without any calculation, that the norm is position independent only when $K$ is restricted always to be 0 which is only true when $L=0$ \ie for the fundamental mode, the general form of which is obtained by linear combination of  components. This completes the proof.

The same conclusion holds for the negative spectrum by interchanging $L$ and $J$.

\begin{ignore}
Thus the positive Reeb field is ${\bde}^M$ and its mirror is the adjoint, ${\bf D}^M$. $M$ is fixed, say $M=0$ and is not now displayed. The notation is that $\big({\bf D}^M\big)_a=\caD^{1M}_a$, and is a function of $q$. Hence construct the combination,
  $$
	     {\bf H}=f{\bde}+g{\bf D}\,,
	$$
where $f$ and $g$ are scalar functions on S$^3$ and require, 
  $$
     \curl {\bf H}=\la\,{\bf H}.
  $$
or, using the rules of the vector calculus,
  $$
	f\curl{\bde}+\grad f\times{\bde}+g\curl{\bf D}+\grad g\times {\bf D}=
	\la(\,f{\bde}+g\,{\bf D}
	$$
or
$$
	\grad f\times{\bde}+\grad g\times {\bf D}=
	(\la-1)(\,f{\bde}+g\,{\bf D}
	$$
which has to be solved for $f$ and $g$	
\end{ignore}
\section{\bf6. Comments and conclusion}

No more examples of the use of spinor spherical harmonics will be given here  as their efficacy has been sufficiently demonstrated by the simple algebraic derivations, presented above, of two theorems on contact structures on the 3--sphere. These are usually obtained by occasionally involved differential geometric means. It is suggested that other, related computations could be similarly eased. Reference [\pref{PSS}] contains a method of constructing particular curl eigenmodes as linear combinations of a Hopf field and its mirror twin which I will leave for another time.

Extensions of the deliberations here could be made to spherical factors, S$^3/\Ga$, when the spectrum ceases to be symmetric, [\pref{Dow11}], as also occurs if the sphere is deformed, [\pref{Dowdewitt}].

I note that related calculations can be found in the topic of knot configurations of rational electromagnetic fields, \eg\ [\pref{KuandL}].

\section{\bf Added Note}

Further reflection has made me modify a comment in this paper, concerning the expansion of vector fields, which I here amplify. 

In works on contact structures on the 3--sphere, a divergence--free vector field is often expanded as $V= \sum_i f_i \al_i$ where $\al_i$ is an orthonormal local frame, the standard choice being the Hopf field and its two partners. These are curl eigenfunctions of eigenvalue $\pm2$ (for a unit 3--sphere). 

If $V$ is a curl eigenfunction of eigenvalue $\mu$ (I choose $\mu>0$)
then it is proved, after some differential analysis, that the $f_i$ are Laplacian eigenfunctions with eigenvalue  $\mu(\mu-2)$. This expansion leads to an associated mapping from S$^3$ to S$^2$ the properties of which, after a slightly involved computation, show, by a contradiction, that, if $\mu>2$, then the norm of $V$ cannot be constant. This constitutes the bulk of the theorem that only the Hopf fields have a constant norm. 
A proof of this has been given given in section 5 which involved only S$^3$ quantities \ie\ without descending to $S^2$.  

Unless I am missing something, it would seem the same conclusion can be reached similarly just from the expansion mentioned above  with the norm squared,
    $
		V^2=\sum_i f_i^2\,
		$
(leading to the above mentioned mapping).		
The scalar functions $f_i$ are constant--coefficient linear combinations of the $\caD^L$ for a fixed $L>0$, if $\mu>2$, so that the $f_i$ are non--constant in this case. Only when $L=0,\,\mu=2$ are they, and therefore $|V|$, constant, which is sufficient.

This direct proof, if it is one, is, in essence, exactly the same as that offered in section 5 only there the expansion is a standard, known construction in angular momentum theory. The advantage of this is that, while a given curl eigenfunction can be written in the form $\sum_i f_i\al_i$, the converse is not true. This renders the {\it construction}, via this form, of useful  curl eigenfunctions somewhat unsystematic, and likewise for the computation of associated Hopf invariants.

\vglue20truept
\noin{\bf References.} \vskip5truept
\begin{putreferences}
  \ref{Dowdewitt}{Dowker,J.S. {\it Vacuum Energy on a Squashed Einstein Universe} in {\it Quantum Gravity}, p.103 edited by  S. C. Christensen (Hilger,Bristol) (1984).}
    \ref{Lyubarskii}{Lyubarskii,G.Ya. {\it The application of Group Theory in Physics}, (Pergamon Press, London) (1960).}
    \ref{KuandL}{Kumar,K. and Lechtenfeld,O. {\it On rational electromagnetic fields} \pla{384}{2020}{126445}.}
     \ref{Wenger}{Wenger,D.L. {\it Representation Functions of the Group of Motions of Clifford Space} \jmp{8}{1967}{135}.}
     \ref{Schr}{Schr\"odinger,E. {\it Maxwell's and Dirac's Equations in the Expanding Universe}, {\it Proc. Roy. Irish Acad.} {\bf A46} (1940) 25.}
		\ref{Schr2}{Schr\"odinger,E. {\it Eigenschwingungen des Sph\"arisches Raumes}, {\it Comm. Pont. Acad. Sci.} {\bf 2} (1938) 321.}
     \ref{Geiges}{Geiges,H.{\it Normal Contact Structures on 3--Manifolds} T$\hat o$hoku Math .J.
		 49 \break (1997) 415.}
     \ref{Pett}{Pettengill,D.F.{\it Spinor Hyperspherical Harmonics and some Applications},
		Ph.D Thesis, Universitiy of Manchester, Manchester (1974).}
    \ref{PSS2}{Peralta--Salas,D. and Slobodeanu,R. {\it Energy minimising Beltrami fields on \break Sasakian 3--manifolds}, {\it Int. Math. Res. Not.} 2021 (2021) 6656; 1806.01164.}
		 \ref{PSS}{Peralta--Salas,D. and Slobodeanu,R. {\it Contact structures and Beltrami fields on the torus and the sphere}, arXiv:2004.10185.}
    \ref{bauer}{Bauer,G. {\it Von den Coefficienten der Reihen von Kugelfunctionen einer Variablen} \jram{56}{1859}{101}.}
    \ref{Biedenharn}{Biedenharn,L.C. {\it An Identity satisfied by Racah coefficients} 
		\jmp{31}{1953}{221}.}
    \ref{niki}{Nikiforov,A.F.,Suslov,S,K. and Uvarov,V.B. {\it Classical Orhogonal Polynomials of a Discrete Variable} (Springer, Berlin (1991)).} 
   \ref{Rayleigh2}{Strutt,J.W. {\it Investigation of the Disturbance produced by a Spherical Obstacle on the Waves of Sound} \plms{4}{1871}{253}.}
   \ref{Rayleigh}{Strutt,J.W. {\it The Theory of Sound}, 2nd. Edn. Vol.II (MacMillan, London (1896)).}
    \ref{Judd} {Judd,B.R. {\it Operator Techniques in Atomic Spectroscopy} (McGraw-Hill, 
		N.Y. (1963)).}
   \ref{FandG}{Freeden,W. and Gutting,M. {\it Special Functions of Mathematical (Geo)-Physics} (Springer, Basel (2013)).}
   \ref{EMOT2}{Erdelyi, A., Magnus, W., Oberhettinger, F. and Tricomi, F.G. {
  \it Higher Transcendental Functions} Vol.2 (McGraw-Hill, N.Y. 1953).}
   \ref{Edmonds}{Edmonds,A.R. {\it Angular Momentum in Quantum Mechanics} (Princeton Univ. Press, Princeton (1957).}
	\ref{Edmonds2}{Edmonds,A.R. {\it Angular Momentum in Quantum Mechanics} (CERN 55-26, Geneva (1955)).}
	\ref{BandT2}{Brussaard,P.J.  and Tolhoek,H.A. {\it Classical L:imits of Clebsch--Gordan Coefficients Racah Coefficients and $D^l_{mn}(\phi,\th,\psi)$--Functions}., {\it Physica} {\bf 23} (1957) 955.}
	 \ref{BandT}{Brussaard,P.J.  and Tolhoek,H.A. {\it Classical L:imits of Clebsch--Gordan Coefficients, Racah Coefficients and $D^l_{mn}(\phi,\th,\psi)$--Functions}., {\it Physica} {\bf 23} (1957) 955.}
	  \ref{FandH}{Frenkel,A. and Hartnoll,S.A. {\it Emergent Area Laws from Entangled Matrices} \break 2301.01325.}
   \ref{YandB}{Yutsis,A.P, and Bandzaitis,A.A. {\it The Theory of Angular Momentum in Quantum Mechanics}, (Mintus, Vilnius (1965).}
   \ref{PandR}{Ponzano,G. and Regge,T. in {\it Spectroscopic and Group Theoretical Methods in Physics}
	(North Holland, Amsterdam (1968).} 
	\ref{Dowk2}{Dowker,J.S. {\it Propagators for Arbitrary Spin in an Einstein Universe}, \aop{71}{1972}{577}.}
    \ref{Nomura}{Nomura,M. {\it Description of the 6--j and 9--j symbols in terms of small numbers of 3--j symbols}. {\it J.Phys.Soc.Japan} {\bf 58} (1989) 2677.}
   \ref{ABHMW}{Alder,K., Bohr,A., Huus,T., Mottelson,B. and Winther,A. {\it Study of Nuclear Structure by Electromagnetic Excitation with Accelerated Ions}, \rmp{28}{1956}{432}.}
   \ref{Dowk2}{Dowker,J.S. {\it Propagators for Arbitrary Spin in an Einstein Universe}, \aop{71}{1972}{577}.}
  \ref{RPS} {Rojo,M.E, Proch\'azka,T. and Sachs,I. {\it On deformations and extensions of Diff(S^2),\break  2105.13375.}}
    \ref{AandM}{Abraham,R. and Marsden,J.E. {\it Foundations of Mechanics}, (Reading, Benjamin /Cummings,1978).}
    \ref{Kowalewski}{Kowalewski,G. {\it Einf\"uhrung in die Theorie der kontinuerlichen Gruppen},
		(Akad. Verlag, Leipzig, 1931).}
    \ref{FandL}{Fradkin,E.S. and Linetsky,V.Ya. {\it Infinite--dimensional generalizations of finite--dimensional symmetries}, \jmp{32}{1991}{1218}.}
    \ref{FandK}{Freidel,L. and Krasnov, K. {\it The fuzzy sphere star product and spin networks}, \jmp{43}{2002}{1737}.}
		\ref{PRS}{Pope,C.N., Romans,L.J. and Shen,X. {\it $W_\infty$  and the Racah--Wigner algebra},\break \np{254}{1991}{401}.}
    \ref{Silberman}{Silberman,L. {\it J.Met} {\bf 11} (1954) 27.}
    \ref{Elsasser}{Elsasser,W.M. {\it Induction Effects in Terrestrial Magnetism Part I. Theory}, \pr{69}{1946}{106}.}
		\ref{ArandS}{Arakelyan,T.A. and Saviddy,G.K. {\it Geometry of a group of area--preserving diffeomorphisms}, \plb{223}{1989}{41}.}
    \ref{Jones}{Jones,M.N. {\it Atmospheric oscillations: I}, {\it Planet. Space Science} {\bf18} (1970) 1393.}
		\ref{Thiebaux}{Thiebaux,M.L. {\it On the Structure of Interaction coefficients in the Spectral Equations for Planetary Waves}, {\it J. Atmospheric Sciences} {\bf28} (1971) 1294.}
    \ref{James1}{James,R.W. {\it The Elsasser and dynamo integrals}, \prs{331}{1973}{469}.}
		 \ref{James2}{James,R.W. {\it The Spectral Form of the Magnetic Induction Equation}, \prs{340}{1974}{287}.}
    \ref{Yoshida}{Yoshida,K. {\it Riemannian curvature on the group of area preserving diffeomorphisms (motion of fluid) of 2-sphere}, {\it Physica} D {\bf 100} (1997) 377.}
		\ref{Moses}{Moses,H.E. {\it The Use of Vector Spherical Harmonics in Global Meteorology and Astronomy}, {\it J. Atmospheric Sciences} {\bf 31} (1974) 1490.}
    \ref{DFMS}{Donnelly,W., Freidel,L., Moosavian,S.F. and Speranza,A.J. {\it Matrix Quantization of Gravitational Edge Modes}, 2212.09120.}
   \ref{Arnold}{Arnol'd, V. {\it Sur la g\'eometrie diff\'erentielle des groupes de Lie de dimension infinie et ses applications \`a l'hydrodynamique des fluides parfaits},  {\it Ann. Inst. Fourier } {\bf 16} (1966) 319.}
	\ref{Dowk3}{Dowker,J.S. {\it Volume preserving diffeomorphisms on the 3--sphere} \cqg{7}{1990}{1241}.}
  \ref{Dowk5}{Dowker,J.S. {\it Diffeomorphisms of the 3--sphere} \cqg{7}{1990}{2353}.}
  \ref{Fock}{Fock,V. \zfp{98}{1935}{145}.}
  \ref{Happer}{Happer,W. \aop{48}{1968}{579}.}
  \ref{FandU}{Falkoff,D.L. and Uhlenbeck,G.E. \pr{79}{1950}{323}.}
  \ref{Rose}{Rose,M.E. \jmp{9}{1962}{409}.}
   \ref{Rose2}{Rose,M.E. \pr{108}{1957}{362}.}
  \ref{Levy}{Levy,M. \prs {204}{1950}{145}.}
  \ref{Schwinger2}{Schwinger,J. \jmp{5}{1964}{1606}.}
  \ref{Muller}{M\"uller,C. {\it Spherical Harmonics} \lnm{17}{1966}{}.}
  \ref{VMK}{Varshalovich,D.A.,Moskalev,A. and Khersonskii,V.K. {\it Quantum Theory of Angular Momentum}, (World Scientific, Singapore, (1988).}
  \ref{DandWo}{Dowker,J.S. and Wolski, A. \prA{46}{1992}{6417}.}
  \ref{Zeitlin1}{Zeitlin,V. {\it Physica D} {\bf 49} (1991).  }
  \ref{Zeitlin0}{Zeitlin,V. {\it Nonlinear World} Ed by
   V.Baryakhtar {\it et al},  Vol.I p.717,  (World Scientific, Singapore, 1989).}
  \ref{Zeitlin2}{Zeitlin,V. \prl{93}{2004}{264501}. }
  \ref{Zeitlin3}{Zeitlin,V. \pla{339}{2005}{316}. }
  \ref{Groenewold}{Groenewold, H.J. {\it Physica} {\bf 12} (1946) 405.}
  \ref{Cohen}{Cohen, L. \jmp{7}{1966}{781}.}
  \ref{AandW}{Argawal G.S. and Wolf, E. \prD{2}{1970}{2161,2187,2206}.}
  \ref{Jantzen}{Jantzen,R.T. \jmp{19}{1978}{1163}.}
  \ref{Moses2}{Moses,H.E. \aop{42}{1967}{343}.}
  \ref{Carmeli}{Carmeli,M. \jmp{9}{1968}{1987}.}
  \ref{SHS}{Siemans,M., Hancock,J. and Siminovitch,D. {\it Solid State
  Nuclear Magnetic Resonance} {\bf 31}(2007)35.}
 \ref{Dowk}{Dowker,J.S. {\it Arbitrary spin theory on the Einstein universe}, \prD{28}{1983}{3013}.}
 \ref{Heine}{Heine, E. {\it Handbuch der Kugelfunctionen}
  (G.Reimer, Berlin. 1878, 1881).}
  \ref{Pockels}{Pockels, F. {\it \"Uber die Differentialgleichung $\De
  u+k^2u=0$} (Teubner, Leipzig. 1891).}
  \ref{Hamermesh}{Hamermesh, M., {\it Group Theory} (Addison--Wesley,
  Reading. 1962).}
  \ref{Racah}{Racah, G. {\it Group Theory and Spectroscopy}
  (Princeton Lecture Notes, 1951). }
  \ref{Gourdin}{Gourdin, M. {\it Basics of Lie Groups} (Editions
  Fronti\'eres, Gif sur Yvette. 1982.)}
  \ref{Clifford}{Clifford, W.K. \plms{2}{1866}{116}.}
  \ref{Story2}{Story, W.E. \plms{23}{1892}{265}.}
  \ref{Story}{Story, W.E. \ma{41}{1893}{469}.}
  \ref{Poole}{Poole, E.G.C. \plms{33}{1932}{435}.}
  \ref{Dickson}{Dickson, L.E. {\it Algebraic Invariants} (Wiley, N.Y.
  1915).}
  \ref{Dickson2}{Dickson, L.E. {\it Modern Algebraic Theories}
  (Sanborn and Co., Boston. 1926).}
  \ref{Hilbert2}{Hilbert, D. {\it Theory of algebraic invariants} (C.U.P.,
  Cambridge. 1993).}
  \ref{Olver}{Olver, P.J. {\it Classical Invariant Theory} (C.U.P., Cambridge.
  1999.)}
  \ref{AST}{A\v{s}erova, R.M., Smirnov, J.F. and Tolsto\v{i}, V.N. {\it
  Teoret. Mat. Fyz.} {\bf 8} (1971) 255.}
  \ref{AandS}{A\v{s}erova, R.M., Smirnov, J.F. \np{4}{1968}{399}.}
  \ref{Shapiro}{Shapiro, J. \jmp{6}{1965}{1680}.}
  \ref{Shapiro2}{Shapiro, J.Y. \jmp{14}{1973}{1262}.}
  \ref{NandS}{Noz, M.E. and Shapiro, J.Y. \np{51}{1973}{309}.}
  \ref{Cayley2}{Cayley, A. {\it Phil. Trans. Roy. Soc. Lond.}
  {\bf 144} (1854) 244.}
  \ref{Cayley3}{Cayley, A. {\it Phil. Trans. Roy. Soc. Lond.}
  {\bf 146} (1856) 101.}
  \ref{Wigner}{Wigner, E.P. {\it Gruppentheorie} (Vieweg, Braunschweig. 1931).}
  \ref{Sharp}{Sharp, R.T. \ajop{28}{1960}{116}.}
  \ref{Laporte}{Laporte, O. {\it Z. f. Naturf.} {\bf 3a} (1948) 447.}
  \ref{Lowdin}{L\"owdin, P-O. \rmp{36}{1964}{966}.}
  \ref{Ansari}{Ansari, S.M.R. {\it Fort. d. Phys.} {\bf 15} (1967) 707.}
  \ref{SSJR}{Samal, P.K., Saha, R., Jain, P. and Ralston, J.P. {\it
  Testing Isotropy of Cosmic Microwave Background Radiation},
  astro-ph/0708.2816.}
  \ref{Lachieze}{Lachi\'eze-Rey, M. {\it Harmonic projection and
  multipole Vectors}. astro- \break ph/0409081.}
  \ref{CHS}{Copi, C.J., Huterer, D. and Starkman, G.D.
  \prD{70}{2003}{043515}.}
  \ref{Jaric}{Jari\'c, J.P. {\it Int. J. Eng. Sci.} {\bf 41} (2003) 2123.}
  \ref{RandD}{Roche, J.A. and Dowker, J.S. \jpa{1}{1968}{527}.}
  \ref{KandW}{Katz, G. and Weeks, J.R. \prD{70}{2004}{063527}.}
  \ref{Waerden}{van der Waerden, B.L. {\it Die Gruppen-theoretische
  Methode in der Quantenmechanik} (Springer, Berlin. 1932).}
  \ref{EMOT}{Erdelyi, A., Magnus, W., Oberhettinger, F. and Tricomi, F.G. {
  \it Higher Transcendental Functions} Vol.1 (McGraw-Hill, N.Y. 1953).}
  \ref{Dowzilch}{Dowker, J.S. {\it Proc. Phys. Soc.} {\bf 91} (1967) 28.}
  \ref{DandD}{Dowker, J.S. and Dowker, Y.P. {\it Proc. Phys. Soc.}
  {\bf 87} (1966) 65.}
  \ref{DandD2}{Dowker, J.S. and Dowker, Y.P. {\it Interactions of Massless Particles of Arbitrary Spin}, \prs{294}{1966}{175}.}
  \ref{CoandH}{Courant, R. and Hilbert, D. {\it Methoden der
  Mathematischen Physik} vol.1 \break (Springer, Berlin. 1931).}
  \ref{Applequist}{Applequist, J. \jpa{22}{1989}{4303}.}
  \ref{Torruella}{Torruella, \jmp{16}{1975}{1637}.}
  \ref{Weinberg}{Weinberg, S.W. \pr{133}{1964}{B1318}.}
  \ref{Meyerw}{Meyer, W.F. {\it Apolarit\"at und rationale Curven}
  (Fues, T\"ubingen. 1883.) }
  \ref{Ostrowski}{Ostrowski, A. {\it Jahrsb. Deutsch. Math. Verein.} {\bf
  33} (1923) 245.}
  \ref{Kramers}{Kramers, H.A. {\it Grundlagen der Quantenmechanik}, (Akad.
  Verlag., Leipzig, 1938).}
  \ref{ZandZ}{Zou, W.-N. and Zheng, Q.-S. \prs{459}{2003}{527}.}
  \ref{Weeks1}{Weeks, J.R. {\it Maxwell's multipole vectors
  and the CMB}.  astro-ph/0412231.}
  \ref{Corson}{Corson, E.M. {\it Tensors, Spinors and Relativistic Wave
  Equations} (Blackie, London. 1950).}
  \ref{Rosanes}{Rosanes, J. \jram{76}{1873}{312}.}
  \ref{Salmon}{Salmon, G. {\it Lessons Introductory to the Modern Higher
  Algebra} 3rd. edn. \break (Hodges,  Dublin. 1876.)}
  \ref{Milnew}{Milne, W.P. {\it Homogeneous Coordinates} (Arnold. London. 1910).}
  \ref{Niven}{Niven, W.D. {\it Phil. Trans. Roy. Soc.} {\bf 170} (1879) 393.}
  \ref{Scott}{Scott, C.A. {\it An Introductory Account of
  Certain Modern Ideas and Methods in Plane Analytical Geometry,}
  (MacMillan, N.Y. 1896).}
  \ref{Bargmann}{Bargmann, V. \rmp{34}{1962}{300}.}
  \ref{Maxwell}{Maxwell, J.C. {\it A Treatise on Electricity and
  Magnetism} 2nd. edn. (Clarendon Press, Oxford. 1882).}
  \ref{BandL}{Biedenharn, L.C. and Louck, J.D. {\it Angular Momentum in Quantum Physics}
  (Addison-Wesley, Reading. 1981).}
  \ref{Weylqm}{Weyl, H. {\it The Theory of Groups and Quantum Mechanics}
  (Methuen, London. 1931).}
  \ref{Robson}{Robson, A. {\it An Introduction to Analytical Geometry} Vol I
  (C.U.P., Cambridge. 1940.)}
  \ref{Sommerville}{Sommerville, D.M.Y. {\it Analytical Conics} 3rd. edn.
   (Bell. London. 1933).}
  \ref{Coolidge}{Coolidge, J.L. {\it A Treatise on Algebraic Plane Curves}
  (Clarendon Press, Oxford. 1931).}
  \ref{SandK}{Semple, G. and Kneebone. G.T. {\it Algebraic Projective
  Geometry} (Clarendon Press, Oxford. 1952).}
  \ref{AandC}{Abdesselam A., and Chipalkatti, J. {\it The Higher
  Transvectants are redundant}, arXiv:0801.1533 [math.AG] 2008.}
  \ref{Elliott}{Elliott, E.B. {\it The Algebra of Quantics} 2nd edn.
  (Clarendon Press, Oxford. 1913).}
  \ref{Elliott2}{Elliott, E.B. \qjpam{48}{1917}{372}.}
  \ref{Howe}{Howe, R. \tams{313}{1989}{539}.}
  \ref{Clebsch}{Clebsch, A. \jram{60}{1862}{343}.}
  \ref{Prasad}{Prasad, G. \ma{72}{1912}{136}.}
  \ref{Dougall}{Dougall, J. \pems{32}{1913}{30}.}
  \ref{Penrose}{Penrose, R. \aop{10}{1960}{171}.}
  \ref{Penrose2}{Penrose, R. \prs{273}{1965}{171}.}
  \ref{Burnside}{Burnside, W.S. \qjm{10}{1870}{211}. }
  \ref{Lindemann}{Lindemann, F. \ma{23} {1884}{111}.}
  \ref{Backus}{Backus, G. {\it Rev. Geophys. Space Phys.} {\bf 8} (1970) 633.}
  \ref{Baerheim}{Baerheim, R. {\it Q.J. Mech. appl. Math.} {\bf 51} (1998) 73.}
  \ref{Lense}{Lense, J. {\it Kugelfunktionen} (Akad.Verlag, Leipzig. 1950).}
  \ref{Littlewood}{Littlewood, D.E. \plms{50}{1948}{349}.}
  \ref{Fierz}{Fierz, M. {\it Helv. Phys. Acta} {\bf 12} (1938) 3.}
  \ref{Williams}{Williams, D.N. {\it Lectures in Theoretical Physics} Vol. VII,
  (Univ.Colorado Press, Boulder. 1965).}
  \ref{Dennis}{Dennis, M. \jpa{37}{2004}{9487}.}
  \ref{Pirani}{Pirani, F. {\it Brandeis Lecture Notes on
  General Relativity,} edited by S. Deser and K. Ford. (Brandeis, Mass. 1964).}
  \ref{Sturm}{Sturm, R. \jram{86}{1878}{116}.}
  \ref{Schlesinger}{Schlesinger, O. \ma{22}{1883}{521}.}
  \ref{Askwith}{Askwith, E.H. {\it Analytical Geometry of the Conic
  Sections} (A.\&C. Black, London. 1908).}
  \ref{Todd}{Todd, J.A. {\it Projective and Analytical Geometry}.
  (Pitman, London. 1946).}
  \ref{Glenn}{Glenn. O.E. {\it Theory of Invariants} (Ginn \& Co, N.Y. 1915).}
  \ref{DowkandG}{Dowker, J.S. and Goldstone, M. \prs{303}{1968}{381}.}
  \ref{Turnbull}{Turnbull, H.A. {\it The Theory of Determinants,
  Matrices and Invariants} 3rd. edn. (Dover, N.Y. 1960).}
  \ref{MacMillan}{MacMillan, W.D. {\it The Theory of the Potential}
  (McGraw-Hill, N.Y. (1930).}
   \ref{Hobson}{Hobson, E.W. {\it The Theory of Spherical and Ellipsoidal Harmonics}
   C.U.P., Cambridge. 1931).}
  \ref{Hobson1}{Hobson, E.W. \plms {24}{1892}{55}.}
  \ref{GandY}{Grace, J.H. and Young, A. {\it The Algebra of Invariants}
  (C.U.P., Cambridge, 1903).}
  \ref{FandR}{Fano, U. and Racah, G. {\it Irreducible Tensorial Sets}
  (Academic Press, N.Y. 1959).}
  \ref{TandT}{Thomson, W. and Tait, P.G. {\it Treatise on Natural Philosophy}
  (Clarendon Press, Oxford. 1867).}
  \ref{Brinkman}{Brinkman, H.C. {\it Applications of spinor invariants in
atomic physics}, North Holland, Amsterdam 1956.}
  \ref{Kramers1}{Kramers, H.A. {\it Proc. Roy. Soc. Amst.} {\bf 33} (1930) 953.}
  \ref{DandP2}{Dowker,J.S. and Pettengill,D.F. {\it The quantum mechanics of the ideal asymmetric top with spin} \jpa{7}{1974}{1527}.}
  \ref{Dowk1}{Dowker,J.S. \jpa{}{}{45}.}
  \ref{DandA}{Dowker,J.S. and Apps, J.S. \cqg{15}{1998}{1121}.}
  \ref{Weil}{Weil,A., {\it Elliptic functions according to Eisenstein
  and Kronecker}, Springer, Berlin, 1976.}
  \ref{Ling}{Ling,C-H. {\it SIAM J.Math.Anal.} {\bf5} (1974) 551.}
  \ref{Ling2}{Ling,C-H. {\it J.Math.Anal.Appl.}(1988).}
 \ref{BMO}{Brevik,I., Milton,K.A. and Odintsov, S.D. \aop{302}{2002}{120}.}
 \ref{KandL}{Kutasov,D. and Larsen,F. {\it JHEP} 0101 (2001) 1.}
 \ref{KPS}{Klemm,D., Petkou,A.C. and Siopsis {\it Entropy
 bounds, monoticity properties and scaling in CFT's}. hep-th/0101076.}
 \ref{DandC}{Dowker,J.S. and Critchley,R. \prD{15}{1976}{1484}.}
 \ref{AandD}{Al'taie, M.B. and Dowker, J.S. \prD{18}{1978}{3557}.}
 \ref{Dow1}{Dowker,J.S. \prD{37}{1988}{558}.}
 \ref{Dow30}{Dowker,J.S. \prD{28}{1983}{3013}.}
 \ref{DandK}{Dowker,J.S. and Kennedy,G. \jpa{}{1978}{}.}
 \ref{Dow2}{Dowker,J.S. \cqg{1}{1984}{359}.}
 \ref{DandKi}{Dowker,J.S. and Kirsten, K. {\it Comm. in Anal. and Geom.
 }{\bf7} (1999) 641.}
 \ref{DandKe}{Dowker,J.S. and Kennedy,G.\jpa{11}{1978}{895}.}
 \ref{Gibbons}{Gibbons,G.W. \pl{60A}{1977}{385}.}
 \ref{Cardy}{Cardy,J.L. \np{366}{1991}{403}.}
 \ref{ChandD}{Chang,P. and Dowker,J.S. \np{395}{1993}{407}.}
 \ref{DandC2}{Dowker,J.S. and Critchley,R. \prD{13}{1976}{224}.}
 \ref{Camporesi}{Camporesi,R. \prp{196}{1990}{1}.}
 \ref{BandM}{Brown,L.S. and Maclay,G.J. \pr{184}{1969}{1272}.}
 \ref{CandD}{Candelas,P. and Dowker,J.S. \prD{19}{1979}{2902}.}
 \ref{Unwin1}{Unwin,S.D. Thesis. University of Manchester. 1979.}
 \ref{Unwin2}{Unwin,S.D. \jpa{13}{1980}{313}.}
 \ref{DandB}{Dowker,J.S.and Banach,R. \jpa{11}{1978}{2255}.}
 \ref{Obhukov}{Obhukov,Yu.N. \pl{109B}{1982}{195}.}
 \ref{Kennedy}{Kennedy,G. \prD{23}{1981}{2884}.}
 \ref{CandT}{Copeland,E. and Toms,D.J. \np {255}{1985}{201}.}
 \ref{ELV}{Elizalde,E., Lygren, M. and Vassilevich,
 D.V. \jmp {37}{1996}{3105}.}
 \ref{Malurkar}{Malurkar,S.L. {\it J.Ind.Math.Soc} {\bf16} (1925/26) 130.}
 \ref{Glaisher}{Glaisher,J.W.L. {\it Messenger of Math.} {\bf18}
(1889) 1.} \ref{Anderson}{Anderson,A. \prD{37}{1988}{536}.}
 \ref{CandA}{Cappelli,A. and D'Appollonio,\pl{487B}{2000}{87}.}
 \ref{Wot}{Wotzasek,C. \jpa{23}{1990}{1627}.}
 \ref{RandT}{Ravndal,F. and Tollesen,D. \prD{40}{1989}{4191}.}
 \ref{SandT}{Santos,F.C. and Tort,A.C. \pl{482B}{2000}{323}.}
 \ref{FandO}{Fukushima,K. and Ohta,K. {\it Physica} {\bf A299} (2001) 455.}
 \ref{GandP}{Gibbons,G.W. and Perry,M. \prs{358}{1978}{467}.}
 \ref{Dow4}{Dowker,J.S..}
  \ref{Rad}{Rademacher,H. {\it Topics in analytic number theory,}
Springer-Verlag,  Berlin,1973.}
  \ref{Halphen}{Halphen,G.-H. {\it Trait\'e des Fonctions Elliptiques},
  Vol 1, Gauthier-Villars, Paris, 1886.}
  \ref{CandW}{Cahn,R.S. and Wolf,J.A. {\it Comm.Mat.Helv.} {\bf 51}
  (1976) 1.}
  \ref{Berndt}{Berndt,B.C. \rmjm{7}{1977}{147}.}
  \ref{Hurwitz}{Hurwitz,A. \ma{18}{1881}{528}.}
  \ref{Hurwitz2}{Hurwitz,A. {\it Mathematische Werke} Vol.I. Basel,
  Birkhauser, 1932.}
  \ref{Berndt2}{Berndt,B.C. \jram{303/304}{1978}{332}.}
  \ref{RandA}{Rao,M.B. and Ayyar,M.V. \jims{15}{1923/24}{150}.}
  \ref{Hardy}{Hardy,G.H. \jlms{3}{1928}{238}.}
  \ref{TandM}{Tannery,J. and Molk,J. {\it Fonctions Elliptiques},
   Gauthier-Villars, Paris, 1893--1902.}
  \ref{schwarz}{Schwarz,H.-A. {\it Formeln und
  Lehrs\"atzen zum Gebrauche..},Springer 1893.(The first edition was 1885.)
  The French translation by Henri Pad\'e is {\it Formules et Propositions
  pour L'Emploi...},Gauthier-Villars, Paris, 1894}
  \ref{Hancock}{Hancock,H. {\it Theory of elliptic functions}, Vol I.
   Wiley, New York 1910.}
  \ref{watson}{Watson,G.N. \jlms{3}{1928}{216}.}
  \ref{MandO}{Magnus,W. and Oberhettinger,F. {\it Formeln und S\"atze},
  Springer-Verlag, Berlin 1948.}
  \ref{Klein}{Klein,F. {\it Lectures on the Icosohedron}
  (Methuen, London, 1913).}
  \ref{AandL}{Appell,P. and Lacour,E. {\it Fonctions Elliptiques},
  Gauthier-Villars,
  Paris, 1897.}
  \ref{HandC}{Hurwitz,A. and Courant,C. {\it Allgemeine Funktionentheorie},
  Springer,
  Berlin, 1922.}
  \ref{WandW}{Whittaker,E.T. and Watson,G.N. {\it Modern analysis},
  Cambridge 1927.}
  \ref{SandC}{Selberg,A. and Chowla,S. \jram{227}{1967}{86}. }
  \ref{zucker}{Zucker,I.J. {\it Math.Proc.Camb.Phil.Soc} {\bf 82 }(1977)
  111.}
  \ref{glasser}{Glasser,M.L. {\it Maths.of Comp.} {\bf 25} (1971) 533.}
  \ref{GandW}{Glasser, M.L. and Wood,V.E. {\it Maths of Comp.} {\bf 25}
  (1971)
  535.}
  \ref{greenhill}{Greenhill,A,G. {\it The Applications of Elliptic
  Functions}, MacMillan, London, 1892.}
  \ref{Weierstrass}{Weierstrass,K. {\it J.f.Mathematik (Crelle)}
{\bf 52} (1856) 346.}
  \ref{Weierstrass2}{Weierstrass,K. {\it Mathematische Werke} Vol.I,p.1,
  Mayer u. M\"uller, Berlin, 1894.}
  \ref{Fricke}{Fricke,R. {\it Die Elliptische Funktionen und Ihre Anwendungen},
    Teubner, Leipzig. 1915, 1922.}
  \ref{Konig}{K\"onigsberger,L. {\it Vorlesungen \"uber die Theorie der
 Elliptischen Funktionen},  \break Teubner, Leipzig, 1874.}
  \ref{Milne}{Milne,S.C. {\it The Ramanujan Journal} {\bf 6} (2002) 7-149.}
  \ref{Schlomilch}{Schl\"omilch,O. {\it Ber. Verh. K. Sachs. Gesell. Wiss.
  Leipzig}  {\bf 29} (1877) 101-105; {\it Compendium der h\"oheren
  Analysis}, Bd.II, 3rd Edn, Vieweg, Brunswick, 1878.}
  \ref{BandB}{Briot,C. and Bouquet,C. {\it Th\`eorie des Fonctions
  Elliptiques}, Gauthier-Villars, Paris, 1875.}
  \ref{Dumont}{Dumont,D. \aim {41}{1981}{1}.}
  \ref{Andre}{Andr\'e,D. {\it Ann.\'Ecole Normale Superior} {\bf 6} (1877)
  265;
  {\it J.Math.Pures et Appl.} {\bf 5} (1878) 31.}
  \ref{Raman}{Ramanujan,S. {\it Trans.Camb.Phil.Soc.} {\bf 22} (1916) 159;
 {\it Collected Papers}, Cambridge, 1927}
  \ref{Weber}{Weber,H.M. {\it Lehrbuch der Algebra} Bd.III, Vieweg,
  Brunswick 190  3.}
  \ref{Weber2}{Weber,H.M. {\it Elliptische Funktionen und algebraische
  Zahlen},
  Vieweg, Brunswick 1891.}
  \ref{ZandR}{Zucker,I.J. and Robertson,M.M.
  {\it Math.Proc.Camb.Phil.Soc} {\bf 95 }(1984) 5.}
  \ref{JandZ1}{Joyce,G.S. and Zucker,I.J.
  {\it Math.Proc.Camb.Phil.Soc} {\bf 109 }(1991) 257.}
  \ref{JandZ2}{Zucker,I.J. and Joyce.G.S.
  {\it Math.Proc.Camb.Phil.Soc} {\bf 131 }(2001) 309.}
  \ref{zucker2}{Zucker,I.J. {\it SIAM J.Math.Anal.} {\bf 10} (1979) 192,}
  \ref{BandZ}{Borwein,J.M. and Zucker,I.J. {\it IMA J.Math.Anal.} {\bf 12}
  (1992) 519.}
  \ref{Cox}{Cox,D.A. {\it Primes of the form $x^2+n\,y^2$}, Wiley,
  New York, 1989.}
  \ref{BandCh}{Berndt,B.C. and Chan,H.H. {\it Mathematika} {\bf42} (1995)
  278.}
  \ref{EandT}{Elizalde,R. and Tort.hep-th/}
  \ref{KandS}{Kiyek,K. and Schmidt,H. {\it Arch.Math.} {\bf 18} (1967) 438.}
  \ref{Oshima}{Oshima,K. \prD{46}{1992}{4765}.}
  \ref{greenhill2}{Greenhill,A.G. \plms{19} {1888} {301}.}
  \ref{Russell}{Russell,R. \plms{19} {1888} {91}.}
  \ref{BandB}{Borwein,J.M. and Borwein,P.B. {\it Pi and the AGM}, Wiley,
  New York, 1998.}
  \ref{Resnikoff}{Resnikoff,H.L. \tams{124}{1966}{334}.}
  \ref{vandp}{Van der Pol, B. {\it Indag.Math.} {\bf18} (1951) 261,272.}
  \ref{Rankin}{Rankin,R.A. {\it Modular forms} C.U.P. Cambridge}
  \ref{Rankin2}{Rankin,R.A. {\it Proc. Roy.Soc. Edin.} {\bf76 A} (1976) 107.}
  \ref{Skoruppa}{Skoruppa,N-P. {\it J.of Number Th.} {\bf43} (1993) 68 .}
  \ref{Down}{Dowker.J.S. \np {104}{2002}{153}.}
  \ref{Eichler}{Eichler,M. \mz {67}{1957}{267}.}
  \ref{Zagier}{Zagier,D. \invm{104}{1991}{449}.}
  \ref{Lang}{Lang,S. {\it Modular Forms}, Springer, Berlin, 1976.}
  \ref{Kosh}{Koshliakov,N.S. {\it Mess.of Math.} {\bf 58} (1928) 1.}
  \ref{BandH}{Bodendiek, R. and Halbritter,U. \amsh{38}{1972}{147}.}
  \ref{Smart}{Smart,L.R., \pgma{14}{1973}{1}.}
  \ref{Grosswald}{Grosswald,E. {\it Acta. Arith.} {\bf 21} (1972) 25.}
  \ref{Kata}{Katayama,K. {\it Acta Arith.} {\bf 22} (1973) 149.}
  \ref{Ogg}{Ogg,A. {\it Modular forms and Dirichlet series} (Benjamin,
  New York,
   1969).}
  \ref{Bol}{Bol,G. \amsh{16}{1949}{1}.}
  \ref{Epstein}{Epstein,P. \ma{56}{1903}{615}.}
  \ref{Petersson}{Petersson.}
  \ref{Serre}{Serre,J-P. {\it A Course in Arithmetic}, Springer,
  New York, 1973.}
  \ref{Schoenberg}{Schoenberg,B., {\it Elliptic Modular Functions},
  Springer, Berlin, 1974.}
  \ref{Apostol}{Apostol,T.M. \dmj {17}{1950}{147}.}
  \ref{Ogg2}{Ogg,A. {\it Lecture Notes in Math.} {\bf 320} (1973) 1.}
  \ref{Knopp}{Knopp,M.I. \dmj {45}{1978}{47}.}
  \ref{Knopp2}{Knopp,M.I. \invm {}{1994}{361}.}
  \ref{LandZ}{Lewis,J. and Zagier,D. \aom{153}{2001}{191}.}
  \ref{DandK1}{Dowker,J.S. and Kirsten,K. {\it Elliptic functions and
  temperature inversion symmetry on spheres} hep-th/.}
  \ref{HandK}{Husseini and Knopp.}
  \ref{Kober}{Kober,H. \mz{39}{1934-5}{609}.}
  \ref{HandL}{Hardy,G.H. and Littlewood, \am{41}{1917}{119}.}
  \ref{Watson}{Watson,G.N. \qjm{2}{1931}{300}.}
  \ref{SandC2}{Chowla,S. and Selberg,A. {\it Proc.Nat.Acad.} {\bf 35}
  (1949) 371.}
  \ref{Landau}{Landau, E. {\it Lehre von der Verteilung der Primzahlen},
  (Teubner, Leipzig, 1909).}
  \ref{Berndt4}{Berndt,B.C. \tams {146}{1969}{323}.}
  \ref{Berndt3}{Berndt,B.C. \tams {}{}{}.}
  \ref{Bochner}{Bochner,S. \aom{53}{1951}{332}.}
  \ref{Weil2}{Weil,A.\ma{168}{1967}{}.}
  \ref{CandN}{Chandrasekharan,K. and Narasimhan,R. \aom{74}{1961}{1}.}
  \ref{Rankin3}{Rankin,R.A. {} {} ().}
  \ref{Berndt6}{Berndt,B.C. {\it Trans.Edin.Math.Soc}.}
  \ref{Elizalde}{Elizalde,E. {\it Ten Physical Applications of Spectral
  Zeta Function Theory}, \break (Springer, Berlin, 1995).}
  \ref{Allen}{Allen,B., Folacci,A. and Gibbons,G.W. \pl{189}{1987}{304}.}
  \ref{Krazer}{Krazer}
  \ref{Elizalde3}{Elizalde,E. {\it J.Comp.and Appl. Math.} {\bf 118}
  (2000) 125.}
  \ref{Elizalde2}{Elizalde,E., Odintsov.S.D, Romeo, A. and Bytsenko,
  A.A and
  Zerbini,S.
  {\it Zeta function regularisation}, (World Scientific, Singapore,
  1994).}
  \ref{Eisenstein}{Eisenstein}
  \ref{Hecke}{Hecke,E.  \ma{112}{1936}{664}.}
  \ref{Hecke2}{Hecke,E. {\it lJber orthogonal-invariante Integralgleichungen} \ma{112} {1918}{398}.}
  \ref{Terras}{Terras,A. {\it Harmonic analysis on Symmetric Spaces} (Springer,
  New York, 1985).}
  \ref{BandG}{Bateman,P.T. and Grosswald,E. {\it Acta Arith.} {\bf 9}
  (1964) 365.}
  \ref{Deuring}{Deuring,M. \aom{38}{1937}{585}.}
  \ref{Guinand}{Guinand.}
  \ref{Guinand2}{Guinand.}
  \ref{Minak}{Minakshisundaram.}
  \ref{Mordell}{Mordell,J. \prs{}{}{}.}
  \ref{GandZ}{Glasser,M.L. and Zucker, {}.}
  \ref{Landau2}{Landau,E. \jram{}{1903}{64}.}
  \ref{Kirsten1}{Kirsten,K. \jmp{35}{1994}{459}.}
  \ref{Sommer}{Sommer,J. {\it Vorlesungen \"uber Zahlentheorie}
  (1907,Teubner,Leipzig).
  French edition 1913 .}
  \ref{Reid}{Reid,L.W. {\it Theory of Algebraic Numbers},
  (1910,MacMillan,New York).}
  \ref{Milnor}{Milnor, J. {\it Is the Universe simply--connected?},
  IAS, Princeton, 1978.}
  \ref{Milnor2}{Milnor, J. \ajm{79}{1957}{623}.}
  \ref{Opechowski}{Opechowski,W. {\it Physica} {\bf 7} (1940) 552.}
  \ref{Bethe}{Bethe, H.A. \zfp{3}{1929}{133}.}
  \ref{LandL}{Landau, L.D. and Lishitz, E.M. {\it Quantum
  Mechanics} (Pergamon Press, London, 1958).}
  \ref{GPR}{Gibbons, G.W., Pope, C. and R\"omer, H., \np{157}{1979}{377}.}
  \ref{Jadhav}{Jadhav,S.P. PhD Thesis, University of Manchester 1990.}
  \ref{DandJ}{Dowker,J.S. and Jadhav, S. \prD{39}{1989}{1196}.}
  \ref{CandM}{Coxeter, H.S.M. and Moser, W.O.J. {\it Generators and
  relations of finite groups} Springer. Berlin. 1957.}
  \ref{Coxeter2}{Coxeter, H.S.M. {\it Regular Complex Polytopes},
   (Cambridge University Press,
  Cambridge, 1975).}
  \ref{Coxeter}{Coxeter, H.S.M. {\it Regular Polytopes}.}
  \ref{Stiefel}{Stiefel, E., J.Research NBS {\bf 48} (1952) 424.}
  \ref{BandS}{Brink, D.M. and Satchler, G.R. {\it Angular momentum theory}, 3rd Edn.
  (Clarendon Press, Oxford, (1993).}
  \ref{Racah3}{Racah G. {\it Theory of Complex Spectra. I}, \pr{61}{1942}{186}.}
  \ref{Schwinger}{Schwinger, J. {\it On Angular Momentum} in {\it Quantum Theory of
  Angular Momentum} edited by Biedenharn,L.C. and van Dam, H.
  (Academic Press, N.Y. (1965)).}
  \ref{Bromwich}{Bromwich, T.J.I'A. {\it Infinite Series},
  (Macmillan, 1947).}
  \ref{Ray}{Ray,D.B. \aim{4}{1970}{109}.}
  \ref{Ikeda}{Ikeda,A. {\it Kodai Math.J.} {\bf 18} (1995) 57.}
  \ref{Kennedy}{Kennedy,G. \prD{23}{1981}{2884}.}
  \ref{Ellis}{Ellis,G.F.R. {\it General Relativity} {\bf2} (1971) 7.}
  \ref{Dow8}{Dowker,J.S. \cqg{20}{2003}{L105}.}
  \ref{IandY}{Ikeda, A and Yamamoto, Y. \ojm {16}{1979}{447}.}
  \ref{BandI}{Bander,M. and Itzykson,C. \rmp{18}{1966}{2}.}
  \ref{Schulman}{Schulman, L.S. \pr{176}{1968}{1558}.}
  \ref{Bar1}{B\"ar,C. {\it Arch.d.Math.}{\bf 59} (1992) 65.}
  \ref{Bar2}{B\"ar,C. {\it Geom. and Func. Anal.} {\bf 6} (1996) 899.}
  \ref{Vilenkin}{Vilenkin, N.J. {\it Special functions},
  (Am.Math.Soc., Providence, 1968).}
  \ref{Talman}{Talman, J.D. {\it Special functions} (Benjamin,N.Y.,1968).}
  \ref{Miller}{Miller, W. {\it Symmetry groups and their applications}
  (Wiley, N.Y., 1972).}
  \ref{Dow3}{Dowker,J.S. \cmp{162}{1994}{633}.}
  \ref{Cheeger}{Cheeger, J. \jdg {18}{1983}{575}.}
  \ref{Dow6}{Dowker,J.S. \jmp{30}{1989}{770}.}
  \ref{Dow20}{Dowker,J.S. \jmp{35}{1994}{6076}.}
  \ref{Dow21}{Dowker,J.S. {\it Heat kernels and polytopes} in {\it
   Heat Kernel Techniques and Quantum Gravity}, ed. by S.A.Fulling,
   Discourses in Mathematics and its Applications, No.4, Dept.
   Maths., Texas A\&M University, College Station, Texas, 1995.}
  \ref{Dow9}{Dowker,J.S. \jmp{42}{2001}{1501}.}
  \ref{Dow7}{Dowker,J.S. \jpa{25}{1992}{2641}.}
  \ref{Warner}{Warner.N.P. \prs{383}{1982}{379}.}
  \ref{Wolf}{Wolf, J.A. {\it Spaces of constant curvature},
  (McGraw--Hill,N.Y., 1967).}
  \ref{Meyer}{Meyer,B. \cjm{6}{1954}{135}.}
  \ref{BandB}{B\'erard,P. and Besson,G. {\it Ann. Inst. Four.} {\bf 30}
  (1980) 237.}
  \ref{PandM}{Polya,G. and Meyer,B. \cras{228}{1948}{28}.}
  \ref{Springer}{Springer, T.A. Lecture Notes in Math. vol 585 (Springer,
  Berlin,1977).}
  \ref{SeandT}{Threlfall, H. and Seifert, W. \ma{104}{1930}{1}.}
  \ref{Hopf}{Hopf,H. \ma{95}{1925}{313}. }
  \ref{Dow}{Dowker,J.S. \jpa{5}{1972}{936}.}
  \ref{LLL}{Lehoucq,R., Lachi\'eze-Rey,M. and Luminet, J.--P. {\it
  Astron.Astrophys.} {\bf 313} (1996) 339.}
  \ref{LaandL}{Lachi\'eze-Rey,M. and Luminet, J.--P.
  \prp{254}{1995}{135}.}
  \ref{Schwarzschild}{Schwarzschild, K., {\it Vierteljahrschrift der
  Ast.Ges.} {\bf 35} (1900) 337.}
  \ref{Starkman}{Starkman,G.D. \cqg{15}{1998}{2529}.}
  \ref{LWUGL}{Lehoucq,R., Weeks,J.R., Uzan,J.P., Gausman, E. and
  Luminet, J.--P. \cqg{19}{2002}{4683}.}
  \ref{Dow10}{Dowker,J.S. \prD{28}{1983}{3013}.}
  \ref{BandD}{Banach, R. and Dowker, J.S. \jpa{12}{1979}{2527}.}
  \ref{Jadhav2}{Jadhav,S. \prD{43}{1991}{2656}.}
  \ref{Gilkey}{Gilkey,P.B. {\it Invariance theory,the heat equation and
  the Atiyah--Singer Index theorem} (CRC Press, Boca Raton, 1994).}
  \ref{BandY}{Berndt,B.C. and Yeap,B.P. {\it Adv. Appl. Math.}
  {\bf29} (2002) 358.}
  \ref{HandR}{Hanson,A.J. and R\"omer,H. \pl{80B}{1978}{58}.}
  \ref{Hill}{Hill,M.J.M. {\it Trans.Camb.Phil.Soc.} {\bf 13} (1883) 36.}
  \ref{Cayley}{Cayley,A. {\it Quart.Math.J.} {\bf 7} (1866) 304.}
  \ref{Seade}{Seade,J.A. {\it Anal.Inst.Mat.Univ.Nac.Aut\'on
  M\'exico} {\bf 21} (1981) 129.}
  \ref{CM}{Cisneros--Molina,J.L. {\it Geom.Dedicata} {\bf84} (2001)
  \ref{Goette1}{Goette,S. \jram {526} {2000} 181.}
  207.}
  \ref{NandO}{Nash,C. and O'Connor,D--J, \jmp {36}{1995}{1462}.}
  \ref{Dows}{Dowker,J.S. \aop{71}{1972}{577}; Dowker,J.S. and Pettengill,D.F.
  \jpa{7}{1974}{1527}; J.S.Dowker in {\it Quantum Gravity}, edited by
  S. C. Christensen (Hilger,Bristol,1984)}
  \ref{Jadhav2}{Jadhav,S.P. \prD{43}{1991}{2656}.}
  \ref{Dow11}{Dowker,J.S. {\it Spherical Universe Topology and the Casimir Effect}, \cqg{21}{2004}4247; hep-th/0404093.}
  \ref{Dow12}{Dowker,J.S. \cqg{21}{2004}4977.}
  \ref{Dow13}{Dowker,J.S. \jpa{38}{2005}1049.}
  \ref{Zagier}{Zagier,D. \ma{202}{1973}{149}}
  \ref{RandG}{Rademacher, H. and Grosswald,E. {\it Dedekind Sums},
  (Carus, MAA, 1972).}
  \ref{Berndt7}{Berndt,B, \aim{23}{1977}{285}.}
  \ref{HKMM}{Harvey,J.A., Kutasov,D., Martinec,E.J. and Moore,G.
  {\it Localised Tachyons and RG Flows}, hep-th/0111154.}
  \ref{Beck}{Beck,M., {\it Dedekind Cotangent Sums}, {\it Acta Arithmetica}
  {\bf 109} (2003) 109-139 ; math.NT/0112077.}
  \ref{McInnes}{McInnes,B. {\it APS instability and the topology of the brane
  world}, hep-th/0401035.}
  \ref{BHS}{Brevik,I, Herikstad,R. and Skriudalen,S. {\it Entropy Bound for the
  TM Electromagnetic Field in the Half Einstein Universe}; hep-th/0508123.}
  \ref{BandO}{Brevik,I. and Owe,C.  \prD{55}{4689}{1997}.}
  \ref{Kenn}{Kennedy,G. Thesis. University of Manchester 1978.}
  \ref{KandU}{Kennedy,G. and Unwin S. \jpa{12}{L253}{1980}.}
  \ref{BandO1}{Bayin,S.S.and Ozcan,M.
  \prD{48}{2806}{1993}; \prD{49}{5313}{1994}.}
  \ref{Chang}{Chang, P. Thesis. University of Manchester 1993.}
  \ref{Barnesa}{Barnes,E.W. {\it Trans. Camb. Phil. Soc.} {\bf 19} (1903) 374.}
  \ref{Barnesb}{Barnes,E.W. {\it Trans. Camb. Phil. Soc.}
  {\bf 19} (1903) 426.}
  \ref{Stanley1}{Stanley,R.P. \joa {49Hilf}{1977}{134}.}
  \ref{Stanley}{Stanley,R.P. \bams {1}{1979}{475}.}
  \ref{Hurley}{Hurley,A.C. \pcps {47}{1951}{51}.}
  \ref{IandK}{Iwasaki,I. and Katase,K. {\it Proc.Japan Acad. Ser} {\bf A55}
  (1979) 141.}
  \ref{IandT}{Ikeda,A. and Taniguchi,Y. {\it Osaka J. Math.} {\bf 15} (1978)
  515.}
  \ref{GandM}{Gallot,S. and Meyer,D. \jmpa{54}{1975}{259}.}
  \ref{Flatto}{Flatto,L. {\it Enseign. Math.} {\bf 24} (1978) 237.}
  \ref{OandT}{Orlik,P and Terao,H. {\it Arrangements of Hyperplanes},
  Grundlehren der Math. Wiss. {\bf 300}, (Springer--Verlag, 1992).}
  \ref{Shepler}{Shepler,A.V. \joa{220}{1999}{314}.}
  \ref{SandT}{Solomon,L. and Terao,H. \cmh {73}{1998}{237}.}
  \ref{Vass}{Vassilevich, D.V. \plb {348}{1995}39.}
  \ref{Vass2}{Vassilevich, D.V. \jmp {36}{1995}3174.}
  \ref{CandH}{Camporesi,R. and Higuchi,A. {\it J.Geom. and Physics}
  {\bf 15} (1994) 57.}
  \ref{Solomon2}{Solomon,L. \tams{113}{1964}{274}.}
  \ref{Solomon}{Solomon,L. {\it Nagoya Math. J.} {\bf 22} (1963) 57.}
  \ref{Obukhov}{Obukhov,Yu.N. \pl{109B}{1982}{195}.}
  \ref{BGH}{Bernasconi,F., Graf,G.M. and Hasler,D. {\it The heat kernel
  expansion for the electromagnetic field in a cavity}; math-ph/0302035.}
  \ref{Baltes}{Baltes,H.P. \prA {6}{1972}{2252}.}
  \ref{BaandH}{Baltes.H.P and Hilf,E.R. {\it Spectra of Finite Systems}
  (Bibliographisches Institut, Mannheim, 1976).}
  \ref{Ray}{Ray,D.B. \aim{4}{1970}{109}.}
  \ref{Hirzebruch}{Hirzebruch,F. {\it Topological methods in algebraic
  geometry} (Springer-- Verlag,\break  Berlin, 1978). }
  \ref{BBG}{Bla\v{z}i\'c,N., Bokan,N. and Gilkey, P.B. {\it Ind.J.Pure and
  Appl.Math.} {\bf 23} (1992) 103.}
  \ref{WandWi}{Weck,N. and Witsch,K.J. {\it Math.Meth.Appl.Sci.} {\bf 17}
  (1994) 1017.}
  \ref{Norlund}{N\"orlund,N.E. \am{43}{1922}{121}.}
  \ref{Duff}{Duff,G.F.D. \aom{56}{1952}{115}.}
  \ref{DandS}{Duff,G.F.D. and Spencer,D.C. \aom{45}{1951}{128}.}
  \ref{BGM}{Berger, M., Gauduchon, P. and Mazet, E. {\it Lect.Notes.Math.}
  {\bf 194} (1971) 1. }
  \ref{Patodi}{Patodi,V.K. \jdg{5}{1971}{233}.}
  \ref{GandS}{G\"unther,P. and Schimming,R. \jdg{12}{1977}{599}.}
  \ref{MandS}{McKean,H.P. and Singer,I.M. \jdg{1}{1967}{43}.}
  \ref{Conner}{Conner,P.E. {\it Mem.Am.Math.Soc.} {\bf 20} (1956).}
  \ref{Gilkey2}{Gilkey,P.B. \aim {15}{1975}{334}.}
  \ref{MandP}{Moss,I.G. and Poletti,S.J. \plb{333}{1994}{326}.}
  \ref{BKD}{Bordag,M., Kirsten,K. and Dowker,J.S. \cmp{182}{1996}{371}.}
  \ref{RandO}{Rubin,M.A. and Ordonez,C. \jmp{25}{1984}{2888}.}
  \ref{BaandD}{Balian,R. and Duplantier,B. \aop {112}{1978}{165}.}
  \ref{Kennedy2}{Kennedy,G. \aop{138}{1982}{353}.}
  \ref{DandKi2}{Dowker,J.S. and Kirsten, K. {\it Analysis and Appl.}
 {\bf 3} (2005) 45.}
  \ref{Dow40}{Dowker,J.S. {\it p-form spectra and Casimir energy}
  hep-th/0510248.}
  \ref{BandHe}{Br\"uning,J. and Heintze,E. {\it Duke Math.J.} {\bf 51} (1984)
   959.}
  \ref{Dowl}{Dowker,J.S. {\it Functional determinants on M\"obius corners};
    Proceedings, `Quantum field theory under
    the influence of external conditions', 111-121,Leipzig 1995.}
  \ref{Dowqg}{Dowker,J.S. in {\it Quantum Gravity}, edited by
  S. C. Christensen (Hilger, Bristol, 1984).}
  \ref{Dowit}{Dowker,J.S. \jpa{11}{1978}{347}.}
  \ref{Kane}{Kane,R. {\it Reflection Groups and Invariant Theory} (Springer,
  New York, 2001).}
  \ref{Sturmfels}{Sturmfels,B. {\it Algorithms in Invariant Theory}
  (Springer, Vienna, 1993).}
  \ref{Bourbaki}{Bourbaki,N. {\it Groupes et Alg\`ebres de Lie}  Chap.III, IV
  (Hermann, Paris, 1968).}
  \ref{SandTy}{Schwarz,A.S. and Tyupkin, Yu.S. \np{242}{1984}{436}.}
  \ref{Reuter}{Reuter,M. \prD{37}{1988}{1456}.}
  \ref{EGH}{Eguchi,T. Gilkey,P.B. and Hanson,A.J. \prp{66}{1980}{213}.}
  \ref{DandCh}{Dowker,J.S. and Chang,Peter, \prD{46}{1992}{3458}.}
  \ref{APS}{Atiyah M., Patodi and Singer,I.\mpcps{77}{1975}{43}.}
  \ref{Donnelly}{Donnelly.H. {\it Indiana U. Math.J.} {\bf 27} (1978) 889.}
  \ref{Katase}{Katase,K. {\it Proc.Jap.Acad.} {\bf 57} (1981) 233.}
  \ref{Gilkey3}{Gilkey,P.B.\invm{76}{1984}{309}.}
  \ref{Degeratu}{Degeratu.A. {\it Eta--Invariants and Molien Series for
  Unimodular Groups}, Thesis MIT, 2001.}
  \ref{Seeley}{Seeley,R. \ijmp {A\bf18}{2003}{2197}.}
  \ref{Seeley2}{Seeley,R. .}
  \ref{melrose}{Melrose}
  \ref{berard}{B\'erard,P.}
  \ref{gromes}{Gromes,D.}
  \ref{Ivrii}{Ivrii}
  \ref{DandW}{Douglas,R.G. and Wojciekowski,K.P. \cmp{142}{1991}{139}.}
  \ref{Dai}{Dai,X. \tams{354}{2001}{107}.}
  \ref{Kuznecov}{Kuznecov}
  \ref{DandG}{Duistermaat and Guillemin.}
  \ref{PTL}{Pham The Lai}
\end{putreferences}

\bye